\documentclass[12pt]{amsart}
\usepackage[dvips]{graphicx,color}
\usepackage[dvips]{epsfig}
\usepackage{amsmath,amssymb,latexsym,amscd,xypic}
\DeclareRobustCommand{\SkipTocEntry}[4]{}
\makeatletter
\newcommand\@dotsep{4.5}
\def\@tocline#1#2#3#4#5#6#7{\relax
  \ifnum #1>\c@tocdepth 
  \else
    \par \addpenalty\@secpenalty\addvspace{#2}%
    \begingroup \hyphenpenalty\@M
    \@ifempty{#4}{%
      \@tempdima\csname r@tocindent\number#1\endcsname\relax
    }{%
      \@tempdima#4\relax
    }%
    \parindent\z@ \leftskip#3\relax \advance\leftskip\@tempdima\relax
    \rightskip\@pnumwidth plus1em \parfillskip-\@pnumwidth
    #5\leavevmode\hskip-\@tempdima #6\relax
    \leaders\hbox{$\m@th
      \mkern \@dotsep mu\hbox{.}\mkern \@dotsep mu$}\hfill
    \hbox to\@pnumwidth{\@tocpagenum{#7}}\par
    \nobreak
    \endgroup
  \fi}
\makeatother 

\DeclareFontFamily{OT1}{rsfs}{}
\DeclareFontShape{OT1}{rsfs}{n}{it}{<-> rsfs10}{}
\DeclareMathAlphabet{\curly}{OT1}{rsfs}{n}{it}

\newcommand\C{\mathbb C}

\newcommand\I{\curly I}

\newcommand\m{\mathfrak m}

\renewcommand\O{\mathcal O}
\newcommand\PP{\mathbb P}
\newcommand\FF{\mathbb F}

\newcommand\Q{\mathbb Q}

\newcommand\Z{\mathbb Z}
\newcommand\mmuu{{\stackrel{\to}{\mu}}}
\newcommand\T{\mathbf{T}}

\newcommand{\rt}[1]{\stackrel{#1\,}{\rightarrow}}
\newcommand{\Rt}[1]{\stackrel{#1\,}{\longrightarrow}}
\newcommand\To{\longrightarrow}
\newcommand\into{\hookrightarrow}
\newcommand\Into{\ar@{^{ (}->}[r]}

\renewcommand\_{^{}_}
\newcommand\take{\backslash}
\newcommand\bull{{\scriptscriptstyle\bullet}}
\newcommand\udot{^\bull}

\newcommand\rk{\operatorname{rank}}
\newcommand\tr{\operatorname{tr}}

\newcommand\Ker{\operatorname{Ker}}

\newcommand\id{\operatorname{id}}
\renewcommand\div{\operatorname{div}}
\renewcommand\Im{\operatorname{Im}}
\newcommand\Hom{\operatorname{Hom}}
\renewcommand\hom{\curly H\!om}
\newcommand\Ext{\operatorname{Ext}}
\newcommand\ext{\curly Ext}

\newcommand\Supp{\operatorname{Supp}}

\newcommand\beq[1]{\begin{equation}\label{#1}}
\newcommand\eeq{\end{equation}}
\newcommand\beqa{\begin{eqnarray*}}
\newcommand\eeqa{\end{eqnarray*}}

\makeatletter \@addtoreset{equation}{section} \makeatother

\newtheorem{defn}[equation]{Definition}
\newtheorem{thm}[equation]{Theorem}

\newtheorem{lem}[equation]{Lemma}
\newtheorem{lems}[equation]{Lemma$^*$}

\newtheorem{prop}[equation]{Proposition}
\newtheorem{conj}[equation]{Conjecture}

\title{{Curve counting via stable pairs in the derived category}}
\author{R. Pandharipande and R. P. Thomas}
\date{April 2009}

\begin{document}

\begin{abstract} \noindent
For a nonsingular projective 3-fold $X$,
we define integer invariants 
virtually enumerating pairs $(C,D)$ where $C\subset X$
is an embedded curve and $D\subset C$ is a divisor.
A virtual class is constructed on the associated moduli space
by viewing a pair as an object in the derived category of $X$.
The resulting invariants are conjecturally equivalent, after
universal transformations, to both the Gromov-Witten and
DT theories of $X$. 
For Calabi-Yau 3-folds, the latter equivalence
should be viewed as a wall-crossing formula in the derived category.

Several calculations of the new invariants are carried out.
In the Fano case, the local contributions of nonsingular
embedded curves are found. In the local toric Calabi-Yau
case, a completely new form of the topological vertex is
described. 

The virtual enumeration of pairs is closely
related to the geometry underlying the BPS
state counts of Gopakumar and Vafa. We prove 
that our integrality predictions for Gromov-Witten
invariants agree with the BPS integrality. Conversely,
the BPS geometry imposes strong conditions on the
enumeration of pairs.
\end{abstract}

\maketitle
\tableofcontents

\setcounter{section}{-1}
\section{Introduction}
There are several ways to compactify the space of embedded curves in
a nonsingular projective variety $X$.
The moduli 
of stable maps $\overline{M}(X)$ provides one such  
compactification.{\footnote{
The locus of embedded curves need not be dense in any of the
compactifications we consider here.}}  The 
Gromov-Witten invariants of $X$ are defined via
integration against the virtual class of $\overline{M}(X)$. 
Because of nontrivial automorphisms,
$\overline{M}(X)$ is a Deligne-Mumford stack, and the Gromov-Witten
invariants are rational numbers. 
Underlying these rational numbers should be integers which,
in a regularized sense, count embedded curves.

Curve counting problems take special form for 3-folds 
since the expected dimension of $\overline{M}(X)$ is independent
of the genus of the map domain. For Calabi-Yau 3-folds, the expected
dimension of $\overline{M}(X)$ is always 0.
The study
of the underlying integer invariants in the Calabi-Yau
case appears already in the quintic 3-fold calculations 
of Candelas, de la Ossa, Green, and Parks \cite{COGP} via 
the Aspinwall-Morrison
formula \cite{AM} in genus 0.

A second compactification of the
space of embedded curves is provided by the Hilbert scheme $I(X)$. 
For 3-folds, $I(X)$ carries a virtual class \cite{ThCasson} and yields
invariants via integration. The resulting theory is
\emph{integer} valued. However, since
$1$-dimensional subschemes of $X$ contain 0-dimensional
subschemes which roam over all of $X$,  the invariants do not directly
count curves. In \cite{MNOP1, MNOP2},
a formal {\em reduced} theory is defined by dividing by the generating series
of 0-dimensional invariants.
For many reasons, a direct geometrical approach to the reduced theory would be 
preferable.

A more recent compactification   
of the space of embedded curves
is provided by the work of Honsen \cite{Honsen}. Honsen defines
 a proper algebraic space $H(X)$ parameterizing
 Cohen-Macaulay curves{\footnote{A Cohen-Macaulay
curve is of pure dimension 1 and possibly
nonreduced, but with no embedded points.}} with finite
maps to $X$ which are generically embeddings. There are no automorphisms
of the map and no roaming 0-dimensional subschemes.
However, a finite number of points
of the curve may be identified in the image in $X$. Unfortunately, 
even for 3-folds, Honsen's
space does not appear to carry 
a natural virtual class.

Honsen's space provides a connection to a fourth compactification. 
The push-forward $f_*\O_C$ of the
structure sheaf associated to an element 
$$[f:C\rightarrow X]\in H(X)$$
is a sheaf on $X$.
What distinguishes $f_*\O_C$
from arbitrary sheaves is the canonical section 
$s\in H^0(X, f_* \O_C)$
obtained from $1\in H^0(C,\O_C)$.
The section $s$ 
may have cokernel at a finite number of points where the curve is not embedded.

We are led to consider the moduli space $P(X)$ of \emph{stable pairs}
 $(F,s)$ where $F$ is a sheaf of fixed Hilbert polynomial 
supported in dimension 1
 and $s\in H^0(X,F)$ is a section. The two stability conditions are: 
\begin{enumerate}
\item[(i)]
the sheaf $F$ is {pure}, 
\item[(ii)] the section $\O_X \stackrel{s}{\rightarrow} F$ has 0-dimensional
cokernel.
\end{enumerate}
By definition, {\em purity} (i) means
every nonzero
subsheaf of $F$ has support of dimension 1 \cite{HLShaves}. In particular,
 purity implies  the (scheme theoretic) 
support $C_F$ of $F$ is a Cohen-Macaulay curve. 
A projective moduli space of stable pairs
 can be constructed by  a standard GIT analysis of Quot scheme
quotients \cite{LPPairs1}. The relationship between our
stability (i)-(ii) and GIT stability is discussed in Section \ref{defs}.

The nicest stable pairs are obtained from the data of an
embedded Cohen-Macaulay curve
$$\iota : C \hookrightarrow X$$
together with a Cartier divisor $D\subset C$. The associated stable pair
is 
\begin{equation}\label{gz2}
(\iota_* \O_C(D),s_D)
\end{equation}
 where
$s_D \in H^0(X,\iota_*\O_C(D))$ is the canonical (up to isomorphism) 
section determined by $D$.
We will often use the abbreviated notation $(\O_C(D),s_D)$
to denote the pair \eqref{gz2}.
A characterization of all stable pairs is provided by
Proposition \ref{desc}.

In case $D$ is empty,
we obtain the pair
$(\O_C,1)$ associated to a Cohen-Macaulay subcurve $C$ alone.
Here, the pair data is equivalent to the kernel of the section
$$
\O_X \stackrel{1}{\rightarrow} \O_C,
$$
which is the ideal sheaf $\I_C$ of $C$. In
fact in $D^b(X)$, the bounded derived category of coherent sheaves on $X$,
the \emph{complex}
$\{\O_X \stackrel{1}{\rightarrow}
\O_C\}$ formed from the pair is quasi-isomorphic to $\I_C$.
More general
stable pairs give rise to more complicated 2-term complexes
\begin{equation*}
I\udot=\{\O_X\stackrel{s}{\rightarrow} F\}\ \ \in\ D^b(X)
\end{equation*}
which may be viewed as Cohen-Macaulay curves $C_F$ obtained
from the kernel of $s$ decorated with a finite number of points
obtained from the cokernel. 
The points replace the 0-dimensional subschemes which appear
 in \cite{MNOP1,MNOP2,ThCasson}. However, they are now
constrained to lie on the curve $C$ rather than being
free to wander over all of $X$.

For a 3-fold $X$,
the natural obstruction theory of  
stable pairs $(F,s)$ 
fails to be 2-term and does \emph{not} admit a virtual class. 
However, in Section \ref{defsec}, we show 
that the fixed-determinant
obstruction theory of the complex 
$I\udot$ in the derived category
provides an {\em alternative} obstruction theory for $P(X)$.
The moduli space $P(X)$ may be naturally viewed
as a well-behaved component of the generally ill-behaved
moduli space of complexes in $D^b(X)$ \cite{LieblichModuli}. Indeed, $P(X)$
provides a rare example where a component of the moduli of complexes
is explicitly constructed as a projective variety.{\footnote{Earlier examples
of well-behaved components of the moduli space of complexes in
the derived category can be found in \cite{BridgelandFlops,InabaModuli}.}}
Arguments parallel to \cite{ThCasson} then show
 the alternative obstruction theory  of $P(X)$ {\em does}
admits a
virtual class of the correct dimension.

Integration against 
the virtual class of $P(X)$ provides
a theory of 3-folds which is deformation invariant and integer valued.
We conjecture our
pair theory for a 3-fold $X$ to be equivalent to the Gromov-Witten theory
and equal to the reduced DT theory.
Integrality constraints on the Gromov-Witten theory of 3-folds
have been predicted by the conjectural BPS invariants 
\cite{GV1,GV2,PandDegen}. 
The integrality predicted by the pair theory and
the BPS invariants are entirely equivalent. 

No direct cohomological definition of BPS invariants satisfying
all expected properties
has yet been proposed.{\footnote{For example, the remarkable
cohomological and
motivic proposals of \cite{hetal,Toda} are unlikely to
be deformation invariant and almost certainly do not satisfy the
Gopakumar-Vafa relationship with Gromov-Witten theory.}}  
Amongst other technical issues, problems with semistability
of sheaves on reducible curves are difficult to overcome. Our moduli
problem stabilizes such sheaves by picking a section. 
The stable pairs invariants are
closely related to the heuristic
interpretation of BPS counts given in \cite{KKVSpinning}.
Indeed, we provide a rigorous definition of BPS 
counts via virtual Euler characteristics modulo a conjectured vanishing
in Section \ref{conjs}.

In Section \ref{examples}, we calculate
the contribution
of a nonsingular embedded curve to the
pair theory in complete agreement
with the Gromov-Witten calculations of \cite{PandDegen}.
In Section \ref{vertex}, 
we present the topological vertex for our invariants. 
The 3-legged topological vertex takes a completely new form
from the pairs point of view. The agreement with
Gromov-Witten theory is still conjectural.

\subsection*{Past work and future directions}
A special case of the moduli space of pairs $P(X)$ arises
naturally 
in Diaconescu's work on local curves
\cite{DiacADHM}. He compactifies a rank 2 vector bundle over a curve 
$C$ to a $\PP^2$-bundle and then uses a relative Beilinson transform on
the fibers to map ideal
sheaves of curves (flat over $C$) to certain quiver sheaves. 
An appropriate stability condition for quiver sheaves then translates
into ours for stable pairs on the $\PP^2$-bundle.

The branch maps of Thaddeus and Alexeev-Knutsen \cite{AKBranch} are \emph{reduced}
curves with finite maps to $X$. 
The moduli space of branch maps $B(X)$ is a proper algebraic space, providing yet another compactification of
the space of embedded curves in $X$. As in the case of Honsen's space, even
in dimension 3, the moduli space $B(X)$ 
seems to lack a natural virtual class.

Following the connection of Honsen's space to $P(X)$, we may
push-forward the structure sheaf of a branch map
$$[f:C \rightarrow X]\in B(X).$$
The result $f_*\O_C$ together with the canonical section $s\in H^0(X,f_*\O_C)$
suggests considering pairs on $X$ consisting
of sheaves of possibly higher rank supported on curves 
with sections whose cokernel may be supported in dimension 1. 
While the construction of such moduli spaces
is easily obtained by varying the stability condition on the space of pairs \cite{LPPairs2},
we have been unable to prove the existence of an appropriate virtual class. 
In particular, the approach to the virtual class of $P(X)$
via derived category deformations discussed in Section \ref{defsec} does
not immediately succeed for higher rank pairs. Perhaps some variant
can be pursued.  

Finally, we point out the Gromov-Witten theory of Calabi-Yau
4 and 5-folds is also known conjecturally
to be governed by integers \cite{KlemmP, PZ}.
Finding equivalent integer valued sheaf theories in higher dimensions is
an interesting problem.

\addtocontents{toc}{\SkipTocEntry}
\subsection*{Acknowledgments.} 
We thank A. J. de Jong for explaining to us the work of his student M. Honsen,
and M. Dr\'ezet for pointing out the work of Le Potier.
Conversations with B. Conrad and M. Lieblich about the deformation theory
of complexes were very helpful. We are grateful to
E. Diaconescu, D. Joyce, and G. Moore for 
conversations about wall crossing formulae in derived categories. We
thank
 J. Bryan, S. Katz, A. Klemm, D. Maulik, and A. Okounkov
for conversations about BPS states and the topological vertex, and the referees for a
thorough reading of the paper and some useful suggestions.

R.P. was partially supported by NSF grant DMS-0500187 and a Packard foundation
fellowship. R.T. was partially supported
by a Royal Society University Research Fellowship.
R.T. would like to thank the Leverhulme Trust and Columbia University for a
visit to New York in the spring of 2007 during which much of the
 work was done.

\section{Definitions} \label{defs}
\subsection{Stability}
Let $X$ be a nonsingular projective 3-fold over $\mathbb{C}$
with a fixed polarization $L$. 
As usual, for sheaves $F$ on $X$,
$$F(k)= F\otimes L^k.$$
Let $q\in\Q[k]$ with positive leading coefficient be a stability
parameter.
For $n\in \mathbb{Z}$ 
 and nonzero
$\beta\in H_2(X,\Z)$, let
$P_{n}^q(X,\beta)$ denote the
moduli space 
of semistable pairs
$$
\O_X\Rt{s}F,
$$
where $F$ is a {pure}  sheaf with Hilbert polynomial
$$\chi(F(k))=k\int_\beta c_1(L) +n$$ and 
$s$ is a nonzero section.
The moduli space $P_{n}^q(X,\beta)$ can be constructed
by GIT \cite{LPPairs2}.

The stability condition for the GIT problem of pairs is defined
as follows.
For a sheaf $G$ with support of dimension at most 1, we let $r(G)$ denote the
coefficient of $k$ in $$\chi(G(k))=r(G)k+c.$$
A {\em proper} subsheaf $G \subset F$ is nonzero and not equal to $F$.
Since $F$ is pure, $G$ has 1-dimensional support and therefore $r(G)>0$.
The pair $(F,s)$ is $q$-stable  if,
for every proper subsheaf $G\subset F$,
\begin{eqnarray} \label{slope}
\frac{\chi(G(k))}{r(G)} &<& \frac{\chi(F(k))+q(k)}{r(F)}\,, \quad \quad
k\gg0
\end{eqnarray}
holds, and for every proper subsheaf $G$ through which 
$s$ factors,
 \begin{eqnarray} \label{slope22}
\frac{\chi(G(k))+q(k)}{r(G)} &<& \frac{\chi(F(k))+q(k)}{r(F)}\,, \quad \quad
k\gg0
\end{eqnarray}
holds. The 
$q$-semistability conditions are obtained from \eqref{slope}-\eqref{slope22} 
after replacing $<$ with $\le$.
 
\subsection{Degree 0}
Le Potier's treatment of the moduli of pairs 
$$\O_X \stackrel{s}{\rightarrow} F$$
is undertaken for
arbitrary Hilbert polynomial{\footnote{Le Potier also treats
arbitrary numbers of sections, not just 1.}} of $F$.
If the class $\beta$ is taken to be 0, 
$\chi(F(k))$ is a constant $n$.
Setting the top order
coefficient $r(F)$ to be $n$, we obtain (semi)-stability conditions
 identical to \eqref{slope}-\eqref{slope22} on proper subsheaves
$G\subset F$.

The $q$-stable pairs with $\beta=0$ 
are precisely those obtained from the structure sheaf 
of a length $n$ subscheme  $S\subset X$ with canonical section,
$$\O_X \stackrel{1}{\rightarrow} \O_S.$$
Indeed, condition (i) is always satisfied and condition (ii)
is excluded since the section 1 generates $\O_S$.
The converse is left to the reader.
Hence, $P^q_{n}(X,0)$ is simply the Hilbert scheme $ $Hilb$(X,n)$ of $n$
points.

\subsection{Limits}
Let $\beta \neq 0$.
We are interested in the moduli space of pairs
$P^q_n(X,\beta)$ in the large $q$ limit. 
 In the degree 1 case,
 $$q(k)=Ak+B,$$ 
the limit is achieved for $A$ and $B$ sufficiently large
(for fixed Hilbert polynomial of $F$).
The limit is
always achieved if $q$ to has degree at least
2.

\begin{lem} \label{stab}
For $q$ sufficiently large, as described above, stability and
semistability coincide. A pair $(F,s)$ is limit stable
if and only if
\begin{enumerate}
\item[(i)]
the sheaf $F$ is {pure}, 
\item[(ii)] the section $\O_X \stackrel{s}{\rightarrow} F$ has 0-dimensional
cokernel.
\end{enumerate}
\end{lem}

\begin{proof}
For $q$ sufficiently large the inequality
\eqref{slope} is always strictly satisfied.
After rearranging inequality (\ref{slope22}) for semistability, we obtain
\beq{slope2}
(r(F)-r(G))\big(\chi(F(k))+q(k)\big)\leq r(F)\ \chi((F/G)(k)),
\eeq
which, for $q$ large, shows both that
\begin{equation}\label{zzz2}
r(F)-r(G)=0
\end{equation} 
and that equality can never occur.
Setting $G=\Im(s)$ in \eqref{zzz2}
implies that $s$ has 0-dimensional cokernel.
\end{proof}

We define a {\em stable pair} $(F,s)$ to be limit stable.
Then, $P_n(X,\beta)$ is the moduli space of stable pairs.
Let 
$$C_F= \Supp(F)\subset X$$
be the scheme theoretic support of $F$.
By condition (ii) of Lemma \ref{stab},  $F$
is isomorphic to the structure sheaf of $C_F$
away from finitely many
points, and so has rank 1 on $C_F$.  

\begin{lem} \label{xx2}
For a stable pair $(F,s)$, the support of $\Im(s)$ is $C_F$.
\end{lem}
\begin{proof}
The issue is local on $X$, so we may consider the geometry on
an affine open on which $F$ is a module.
The supports of $F$ and $\Im(s)$   are defined by the annihilators of $F$ and $s$ respectively.
The annihilator of $F$ certainly annihilates $s$.  
Conversely, let $a\in$\,Ann$(s)$ be a function. 
If $a\notin\,$Ann$(F)$, let $f\in F$ be a section for which $af\in F$
does not vanish. 
Then, the submodule of $F$ generated by $af$
has dimension 0 support (away from the nonempty open set
on which $s$ generates $F$ guaranteed by condition (ii)
of pair stability). Hence, the purity of $F$ is violated.
\end{proof}

Since $\Im(s)$ is a quotient of $\O_X$, $\Im(s)$ is a structure
sheaf. By Lemma \ref{xx2}, $\Im(s)\cong\O_{C_F}$.
As a subsheaf of a pure sheaf, $\Im(s)$
 is also pure. Therefore,
 $C_F$ is Cohen-Macaulay.

The following kernel/cokernel exact sequence is associated to the
stable pair $(F,s)$,
\beq{IOFQ}
0\to\I_{C_F}\to\O_X\Rt{s}F\to Q\to0.
\eeq
The cokernel
$Q$ has dimension 0 support by stability.
The {\em reduced} support scheme, $\Supp^{red}(Q)$, is 
called the {\em zero locus} of the pair.
The zero locus lies on $C_F$.

Let $C\subset X$ be a fixed Cohen-Macaulay curve, and $\m\subset\O_C$ the ideal of a
finite union of closed points. We now characterize stable pairs with support $C$ and
zero locus supported at these points.

Since $$\hom(\m^r/\m^{r+1},\O_C)=0$$ by
the purity of $\O_C$, we obtain an inclusion $$\hom(\m^r,\O_C)\subset
\hom(\m^{r+1},\O_C).$$ 
The inclusion $\m^r\into\O_C$ induces a canonical section 
$$\O_C\into\hom(\m^r,\O_C).$$

\begin{prop} \label{desc}
A stable pair $(F,s)$ with support $C$ satisfying
$$\Supp^{red}(Q) \subset \Supp(\O_C/\m)$$ 
is equivalent to a subsheaf
of $\hom(\m^r,\O_C)/\O_C,\ r\gg0.$
\end{prop}

Alternatively, we may work with coherent subsheaves of the quasi-coherent sheaf
$\lim\limits_{\To}\hom(\m^r,\O_C)/\O_C$.

\begin{proof}
Let $Q$ denote the $0$-dimensional cokernel of the stable pair. Its dual $\hom(Q,\O_C)$
vanishes, since $\O_C$ is pure. Therefore, applying $\hom(\
\cdot\ ,\O_C)$ to
\beq{Q}
0\to\O_C\to F\to Q\to0
\eeq
yields the inclusion
\beq{Fdual}
0\to\hom(F,\O_C)\to\O_C.
\eeq
Hence, $\hom(F,\O_C)$ is the pushfoward to $X$ of an ideal sheaf
$\I_Z$ on $C$. 
Since (\ref{Fdual}) is a generic isomorphism, $Z$ is 
0-dimensional and 
$$Z^{red}\subset \Supp^{red}(Q).$$
Dualizing on $C$ again gives
\beq{F}
0\to\O_C\to\hom(\I_Z,\O_C).
\eeq
The obvious double dual map 
$$F\to\hom(\hom(F,\O_C),\O_C)=\hom(\I_Z,\O_C)$$ is generically
an isomorphism, so is an injection by the purity of $F$. The map (\ref{F})
factors through the original section $\O_C\to F$. Thus, we have the data
$$
\O_C\to F\subseteq\hom(\I_Z,\O_C),
$$
with the composition being the canonical section of $\hom(\I_Z,\O_C)$.

For $r\gg0$, there is an inclusion $\m^r\subset\I_Z$ with 0-dimensional cokernel. Therefore by purity we get inclusions
$$\hom(\I_Z,\O_C)\subset\hom(\m^r,\O_C)\subset\hom(\m^{r+1},\O_C).$$
We
obtain the subsheaf
\beq{data}
F\subset\lim_{\To}\hom(\m^r,\O_C),
\eeq
containing the canonical section of $\hom(\m^r,\O_C)$. Dividing by the
canonical section
gives a coherent subsheaf $$Q\subset\lim\limits_{\To}\hom(\m^r,\O_C)/\O_C.$$

Conversely, given such a $Q$ in $\hom(\m^r,\O_C)/\O_C$, the
sheaf $F$ is recovered as its inverse image \eqref{data} in $\hom(\m^r,\O_C)$.
 Moreover, $F$ has a section $s$ fitting into an exact sequence 
\eqref{Q}. As a subsheaf of
$\hom(\m^r,\O_C)$, which is pure since $C$ is Cohen-Macaulay, $F$ is also pure.
The pair $(F,s)$ supported
on $C$ is stable by Lemma \ref{stab}.
\end{proof}

\subsection{Derived category}\label{dcat}

Let $D^b(X)$ be the bounded derived category of coherent sheaves on $X$.
Let $I\udot\in D^b(X)$ be determined by the complex
 $$\{\O_X\Rt{s}F\}$$
associated to the stable pair $(F,s)$
 with $\O_X$ in degree 0.
In what follows, all Hom and Ext groups are considered in $D^b(X)$.

 We have the following exact triangles in $D^b(X)$ associated to $I\udot$:
\begin{eqnarray} \label{a}
F[-1]\to I\udot\to\O_X\rt{s}F\to\ldots, \\
\I_C\to I\udot\to Q[-1]\to\I_C[1]\to\ldots\,, \label{b}
\end{eqnarray}
the second coming from (\ref{IOFQ}).

\begin{lem} \label{homs}
$\ext^{\le-1}(I\udot,I\udot)=0$ and $\hom(I\udot,I\udot)=\O_X$.
\end{lem}

\begin{proof}
Applying $\hom(\ \cdot\ ,\O_X)$ to (\ref{a}) yields
\begin{multline}\label{homo}
\hom(F,\O_X) \to\hom(\O_X,\O_X)\to 
\hom(I\udot,\O_X)\\
\to\ext^1(F,\O_X).
\end{multline}
The first and last terms vanish since $F$ has support of codimension $2$. 
The identity generates the second term, and maps in the third term to 
the canonical map $I\udot\to\O_X$ of \eqref{a}. This canonical map therefore generates
\begin{equation}\label{ccv}
\hom(I\udot,\O_X)\cong\O_X. 
\end{equation}
It is the image of the identity in the exact sequence
\beq{p}
\ext^{-1}(I\udot,F)\to\hom(I\udot,I\udot)\to\hom(I\udot,\O_X)
\eeq
obtained from \eqref{a} by applying $\hom(I\udot,\ \cdot\ )$.

Therefore to show that $\hom(I\udot,I\udot)=\O_X$ we need only prove the vanishing of
$\ext^{-1}(I\udot,F)$.
But $\hom(\ \cdot\ ,F)$ applied to (\ref{b}) gives
\begin{equation*}
\ext^{-1}(I\udot,F)
\cong\hom(Q,F),
\end{equation*}
which vanishes by the purity of $F$.

The same sequences in lower degrees prove the
vanishing of the sheaves $\ext^{\le-1}(I\udot,I\udot)$.
\end{proof}

After tensoring the result of Lemma \ref{homs} by $K_X$, we obtain
\begin{equation} \label{ddd}
 \ext^{\le -1}(I\udot,I\udot\otimes K_X)=0 \text{  and  }
\hom(I\udot,I\udot \otimes K_X)=K_X.
\end{equation}

\begin{lem} \label{dddd}$\ext^{\le -1}(I\udot,\O_X)=0$ and
$\hom(I\udot,\O_X)=\O_X$.
\end{lem}

\begin{proof} 
The second claim is \eqref{ccv}.
The vanishing of $\ext^{\le-1}(I\udot,\O_X)$ is obtained
from \eqref{homo} in lower degrees. 
\end{proof}

By the local-to-global spectral sequence and Lemma \ref{dddd}, we obtain
$${\Hom}(I\udot,\O_X)=\C.$$
The exact triangle \eqref{a} is obtained canonically
from the unique (up to scalars)
nonzero element  $ {\Hom}(I\udot,\O_X)$ with $F$ quasi-isomorphic
to the mapping cone $M$.
The pair $(F,s)$ can be recovered from the complex $I\udot\in D^b(X)$
from the $0^{th}$ cohomology of the induced map
$$\O_X \rightarrow M.$$
Hence, by considering $I\udot \in D^b(X)$, no information about the original
pair is lost. We have proven the following result.

\begin{prop}
The stable pairs $(F,s)$ and $(F',s')$ are isomorphic if and only
if the complexes
$\{\O_X\Rt{s}F\}$ and $\{\O_X\Rt{s'}F'\}$
are quasi-isomorphic.
\end{prop}

\section{Deformation theory and the virtual class} \label{defsec}
\subsection{Pairs and complexes}
Let $X$ be a 3-fold. As before, to each stable pair
\begin{equation*}
[\O_X\stackrel{s}{\to} F]\in P_n(X,\beta)
\end{equation*}
we associate a complex
\begin{equation*}
I\udot=\{\O_X\to F\}\in D^b(X).
\end{equation*}
The first order deformation theory of the moduli space $P_{n}(X,\beta)$ 
of stable pairs is governed
by the tangent space $\Ext^0(I\udot,F)$ and 
the obstruction space $\Ext^1(I\udot,F)$ \cite{LPPairs1}.
The deformation theory of the complex $I\udot$ \emph{with fixed determinant $\O_X$}
in $D^b(X)$ 
is governed by $\Ext^1(I\udot,I\udot)_0$ and $\Ext^2(I\udot,I\udot)_0$ \cite{InabaModuli,
LieblichModuli}. The subscript 0 in the latter two groups denotes trace-free Ext. 

Since a deformation of the pair $(F,s)$
induces a deformation of the complex  $I\udot$
with trivial determinant, there is a map 
from the first pair of
groups to the second,
$$
\Ext^i(I\udot,F)\to\Ext^{i+1}(I\udot,I\udot)_0\,,
$$
obtained by applying $\Hom(I\udot,\ \cdot\ )$ to the canonical map $F[-1]\to I\udot$ of (\ref{a}) and showing the image is in the trace-free part of
$\Ext^*(I\udot,I\udot)$. 

\subsection{Tangent spaces}
We first show that to all orders the deformations of pairs $(F,s)$ equal the deformations of complexes $I\udot$ of fixed determinant.
Hence, $P_n(X,\beta)$ is 
a locally 
complete moduli space of \emph{complexes} $I\udot$ of fixed determinant and
the map
$$
\Ext^0(I\udot,F)\to\Ext^{1}(I\udot,I\udot)_0
$$
is an isomorphism.

A {\em family of stable pairs} over a quasi-projective base scheme $B$ is a pair
$$
\O_{X\times B}\stackrel{s}{\to} F,
$$
on $X\times B$,
for which  
\begin{enumerate}
\item[(i)]
$F$ is flat over $B$, 
\item[(ii)] for all closed points $b \in B$,
the restriction $(F_b,s_b)$ to the 
fiber $X\times\{b\}$ is a stable pair. 
\end{enumerate}
Let $Q$ be the cokernel of $s$. The sheaf $Q$ is supported
in relative dimension 0 over $B$. 

A {\em family of complexes} over a base $B$ is a {\em perfect complex}
$I\udot$ on $X \times B$.
By definition, a perfect complex has a 
finite resolution by locally free sheaves.\footnote{Since $X\times B$ 
is quasi-projective, the local and global existence of a  resolution
are equivalent.}
No additional flatness condition over $B$ is required.

For the deformation question, we study how families extend over 
nilpotent thickenings of the base. 
Let 
$$B\supset B_0$$ 
be a thickening of base schemes
where $B_0$ is defined
by a nilpotent ideal $J$ (satisfying $J^N=0$ for some $N>0$). 
A family of stable pairs $(F,s)$ over $B$ is a deformation of a family 
$(F_0,s_0)$ over $B_0$ if the restriction of $(F,s)$ over $B_0$ is
isomorphic to $(F_0,s_0)$.
A family of complexes $I\udot$ over $B$ is a deformation of a family
$I\udot_0$ over $B_0$ if the
derived restriction of $I\udot$ to $X\times B_0$ is quasi-isomorphic
to $I_0\udot$.

Let $(F_0,s_0)$ be a family of stable pairs over $B_0$.
Let
$$I\udot_0=\{\O_X\to F_0\}$$
be the associated family of complexes over $B_0$. The flatness of $F_0$
implies that $I_0\udot$ is perfect.

\begin{lem} \label{itoiv}
Every deformation $I\udot$ over $B$ of $I\udot_0$
 is quasi-isomorphic to a 2-term complex
of sheaves $\{A\to E^1\}$ satisfying
\begin{enumerate}
\item[(i)] $E^1$ is locally free, 
\item[(ii)]
$A$ is a pure sheaf of local projective dimension at most 1,
\item[(iii)] $h^1(I\udot)$ has support of relative
dimension 0 over $B$,
\item[(iv)]  
$h^0(I\udot)$, away from the support of $h^1(I\udot)$,
 is flat over $B$.
\end{enumerate}
\end{lem}

\begin{proof} Since $F_0$ is pure of relative dimension 1
 on $$\pi_0\colon X\times B_0 \rightarrow B_0,$$ $F_0$ 
has depth 1 on the fibers of $\pi_0$.
By the Auslander-Buchsbaum formula, $F_0$ has local projective
 dimension at most 2 on the fibers of $\pi_0$.
Let $$K^\bullet \rightarrow F_0 \rightarrow 0$$
be a locally free resolution of $F_0$.
Consider the cut-off
$$0\rightarrow \Ker  \rightarrow K^{-1} \rightarrow K^0 \rightarrow F_0 
\rightarrow 0.$$
By flatness of $F_0$ and the projective dimension computation, 
the restriction of $\Ker$ to each fiber $\pi_0^{-1}(b_0)$
is locally free. Since $\Ker$ is also flat over $B_0$, we conclude
$\Ker$ is locally free. Hence,
$F_0$ has a locally free resolution of length 3
on $X\times B_0$. 
Since $\O_{X\times B_0}$ is already locally free, $I\udot_0$ is 
quasi-isomorphic to a complex of locally free sheaves of length 3 on $X$.

Let $E\udot$
be a finite complex of locally free sheaves on $X\times B$
  quasi-isomorphic to $I\udot$.
We will use standard base change, semicontinuity, and Nakayama Lemma 
arguments to trim the complex $E\udot$ down to length 3. 

Let $E^n$ be 
the last nonzero term of $E\udot$. If $n>1$,
then 
$$E^{n-1}|_{X\times B_0}\to E^n|_{X\times B_0}$$ is surjective. 
So 
$
E^{n-1}\to E^n
$
is surjective in a neighbourhood of $X\times B_0$ and thus
on all of $X\times B$. The kernel of
$
E^{n-1}\to E^n
$
is then locally
free, and
$E\udot$ can be trimmed. We can thus assume that $E^1$ is the last term.

Similarly, let $E^m$ be the first nonzero term.
If $m<-1$, then 
$$E^m|_{X\times B_0}\to E^{m+1}|_{X\times B_0}$$ is injective on {\em fibers}
by base change and the 3 term result for $I\udot_0$.
So $$E^m\to E^{m+1}$$ is injective on fibers with locally free cokernel,
and $E\udot$ may be again trimmed. We therefore assume that $E^{-1}$ is the first term.

We conclude that $I\udot$ is quasi-isomorphic to a length 3 complex of locally
free sheaves 
$$E^{-1}\to E^0\to E^1$$
on $X\times B$. The first map
$$E^{-1}|_{X\times B_0}\to E^0|_{X\times B_0}$$ is injective
as a map of sheaves since $h^{-1}(I\udot_0)=0$.
Hence,
 $E^{-1}\to E^0$ is injective in a neighbourhood of $X\times B_0$ and thus
 on all of $X\times B$. 
So $I\udot$ is quasi-isomorphic to $$A\to E^1$$ 
for some sheaf $A=E^0/E^{-1}$ of projective
 dimension at most 1. 
Moreover, $E^{-1}\to E^0$ is actually injective on fibers 
away from the relative curve $$C_0\times B_0\subset
X\times B$$ on which $F_0$ is supported. Hence,
$A$ is locally free away from $C_0$.

To establish the purity of $A$, we must show that any subsheaf
$A'\subset A$ with support of codimension at least $1$ is in fact zero.
Since $A$ is locally free away
from $C_0$, $A'$ must be supported on $C_0$ and so has codimension at 
least $2$. Let $\pi$ denote the projection 
$$\pi\colon X\times B \to B$$
and let $K_X$ be the relative dualizing sheaf.
Then, 
$$\Hom(A',A)=H^0(\pi_*\hom(A',A)),$$ and we need only prove 
$\pi_*\hom(A',A)=0$.

By relative Serre duality for the smooth map $\pi$, the derived dual of
$R\pi_*R\hom(A',A)[3]$ is quasi-isomorphic to
$$
R\pi_*R\hom(R\hom(A',A),K_X)=R\pi_*(R\hom(A,A'\otimes K_X)).
$$
The latter's $k$th cohomology sheaf
can be calculated by the local-to-global spectral sequence
with $E_2$ term
\begin{equation}\label{vv551}
R^i\pi_*\,\ext^j(A,A'\otimes K_X), \qquad i+j=k.
\end{equation}
As $A$ has projective dimension at most $1$, we have vanishing for $j\ge2$. Since
$A'$ is supported in relative dimension at most $1$, we have vanishing
for $i\ge2$.  Therefore, the sheaves \eqref{vv551} vanish for $i+j\ge3$.

It follows that $R\pi_*R\hom(A,A'\otimes K_X)$ is quasi-isomorphic to a complex supported in degrees $0,1$ and $2$. Taking the derived
dual gives a complex in degrees $\ge-2$. Shifting by $[-3]$ then shows that $R\pi_*R\hom(A',A)$
is supported in degrees $\ge1$. Therefore its $0$th degree cohomology sheaf $\pi_*\hom(A',A)$
vanishes, and $A$ is indeed pure.

The complement of the closed subscheme $Z_0=\Supp(Q_0)$ determines
open sets
$$U_0=(X\times B_0)\take Z_0, \ \ U=(X\times B)\take Z_0.$$
Though $U_0$ and $U$ have the same closed points, $U_0$ is a closed
subscheme of $U$. 
Certainly,
$$
E^0|_{U_0}\to E^1|_{U_0} 
$$
is surjective. By 
Nakayama's Lemma, 
\begin{equation}\label{hhttj}
E^0|_{U}\to E^1|_{U} 
\end{equation}
is also surjective.
Thus, $h^1(I\udot)$ has relative dimension 0 support, which is property
(iii).

Let $K|_U$ be the locally free kernel of \eqref{hhttj}. 
The sheaf
$h^0(I\udot)|_{U}$ is quasi-isomorphic to 
$$E^{-1}|_{U}\to K|_U$$
by the established injectivity.
The complex
\begin{equation}\label{bnnm}
E^{-1}|_{U_0}\to K|_{U_0}
\end{equation}
is the derived
restriction of $h^0(I\udot)|_{U}$ to $U_0$.
By repeating the argument for $B_0$ instead of $B$, we
find \eqref{bnnm} is
quasi-isomorphic to its cokernel
 $$h^0(I\udot_0)|_{U_0}\cong h^0(I\udot)|_{U_0}.$$

Since
$h^0(I\udot_0)|_{U_0}$
is the kernel of the surjection 
$$\O_{U_0} \rightarrow F_0|_{U_0}$$
and $F_0$ is flat over $B_0$, we see
that $h^0(I\udot_0)|_{U_0}$ is flat over $B_0$.
By Lemma \ref{flatnesslem} below,
$h^0(I\udot)|_{U}$ is flat over $B$, which is (iv).
\end{proof}

\begin{lem} \label{flatnesslem}
Let $\iota\colon B_0\into B$ be a nilpotent thickening,
 and let
$F$ be a coherent sheaf on an
open set 
$$U\subset X\times B.$$ Let $F_0=\iota^*F$ be the
restriction to $U_0=\iota^*(U)$.
Then, $F$ is flat over $B$ if and only if
\begin{enumerate}
\item[(i)] $L\iota^*F\cong F_0$ and
\item[(ii)] $F_0$ is flat over $B_0$.
\end{enumerate}
\end{lem}

\begin{proof}
Flatness clearly implies (i) and (ii). For the converse we must show that
\beq{ffflat}
F\stackrel{L}{\otimes}M\cong F\otimes M
\eeq
for any $\O_B$-module $M$.
For $M=\iota_*M_0$, where $M_0$ is an $\O_{B_0}$-module, \eqref{ffflat}
 is clear:
$$
F\stackrel{L}{\otimes}\iota_*M_0\cong\iota_*(L\iota^*F\stackrel{L}{\otimes}M_0)
\cong\iota_*(F_0\stackrel{L}{\otimes}M_0)\cong\iota_*(F_0\otimes M_0)
\cong F\otimes\iota_*M_0.
$$
Here, the second isomorphism comes from (i) and the third from (ii).

Since $B$ is a nilpotent thickening of $B_0$, $M$ can be written as a finite
series of extensions of such $\O_{B_0}$-modules, which gives (\ref{ffflat}).
\end{proof}

\begin{thm} \label{def}
Every deformation $I\udot$ over $B$ of $I\udot_0$
with trivial determinant 
is quasi-isomorphic
to a complex $$\{\O_{X\times B}\rt{s}F\},$$
where $F$ is a flat deformation of $F_0$ with section $s$.
\end{thm}

\begin{proof}
The rank of $I\udot$ is 1.
 Since $$Q=h^1(I\udot)$$ has rank 0 on $X\times B$,
the rank of
$h^0(I\udot)$ must be $1$. 
As a subsheaf of $A$, $h^0(I\udot)$ is also pure
so injects into its double dual,
$$0 \rightarrow h^0(I\udot)\rightarrow h^0(I\udot)^{\vee \vee}.$$
By (iii) and (iv) of Lemma \ref{itoiv}, 
$h^0(I\udot)$ is flat over $B$
away from a set of codimension 3.
Therefore by \cite[Lemma 6.13]{KollarJDG} the double dual is locally free
away from the codimension 3 set. Moreover, in the proof of \cite[Lemma
6.13]{KollarJDG}, Koll\'ar shows that a reflexive rank 1 sheaf which is locally
free away from a set of codimension 3 is in fact locally free globally. Finally
$h^0(I\udot)^{\vee \vee}$ has trivial determinant since $I\udot$ does, so
we find that
$$h^0(I\udot)^{\vee \vee} \cong \O_{X\times B}$$
and $h^0(I\udot)$ is an ideal sheaf $\I_C\subset\O_{X\times B}$.

The exact triangle $h^0(I\udot)\to I\udot\to h^1(I\udot)[-1]$, written as
$$
\I_C\to I\udot\to Q[-1],
$$
describes $I\udot$ as the cone of a map $Q[-2]\to\I_C$. The
latter is an element
$\alpha\in\Ext^2(Q,\I_C)$.

Consider the exact sequence obtained from the ideal sequence
\begin{equation}\label{zdg}
0\rightarrow \I_C \rightarrow
\O_{X\times B} \rightarrow \O_C \rightarrow 0
\end{equation}
by applying $\Hom(Q,\ \cdot\ )$:
$$
\Ext^1(Q,\O_{X\times B})\to\Ext^1(Q,\O_C)\to\Ext^2(Q,\I_C)\to
\Ext^2(Q,\O_{X\times B}).
$$
The first and last terms vanish since $Q$ has codimension 3 support. 
Thus, $$\Ext^1(Q,\O_C)\cong\Ext^2(Q,\I_C),$$ 
and $\alpha$ is the cup product of an element
$\epsilon\in\Ext^1(Q,\O_C)$ 
with the extension class in $\Ext^1(\O_C,\I_C)$ of the ideal sequence 
\eqref{zdg}. 
Thus it is represented by the splicing together of the exact sequences
$$
\spreaddiagramcolumns{-1pc}
\spreaddiagramrows{-1.2pc}
\xymatrix{
&&& \qquad 0 \rto & \O_C \ar@{=}[d]\rto & F\rto & Q \rto & 0 \\
\hspace{-2cm} \mathrm{} & 0 \rto & \I_C \rto & \O_{X\times B} \rto & \O_C \rto & 0\,,
}$$
where $F$ is the extension defined by $\epsilon$. 
The result is the exact sequence
$$
0\to\I_C\to\O_{X\times B}\to F\to Q\to0.
$$
Hence, $I\udot$ is quasi-isomorphic to the complex 
$$
\{\O_{X\times B}\to F\}.
$$

Finally, we show $F$ is flat over $B$. 
By Lemma \ref{flatnesslem}, flatness follows if
$$
L\iota^*F\cong F_0,
$$
where $\iota\colon B_0\into B$ is the inclusion.

$F$ can be described as the cone of the canonical map \eqref{a},
 $$I\udot\to\O_{X\times B},$$  so the derived restriction of $F$
to $X\times B_0$ is the cone on the induced map 
$$I_0\udot\to\O_{X\times B_0}.$$ 
Similarly $F_0$ is the cone on such a map \eqref{a}, 
so we need only check the two maps coincide.

By sequence \eqref{homo}, $\hom(I\udot,\O_{X\times B})\cong\O_{X\times B}$.
In lower degrees, \eqref{homo} yields the vanishing  
of
$\ext(I\udot, \O_{X\times B})^{\leq -1}$. Therefore,
 $$\Hom(I\udot,\O_{X\times B})\cong\Gamma(\O_B)$$ 
is generated over $\O_B$ by the canonical map (\ref{a}). Similarly,
$$\Hom(I_0\udot,\O_{X\times B_0})\cong\Gamma(\O_{B_0})$$ 
is generated by the
canonical map.
The restriction to $X\times B_0$ of the canonical map on $X\times B$ is 
$\phi$ times the canonical map on $X\times B_0$ for some $\phi\in
\Gamma(\O_{B_0})$. However, away from the codimension 2 support $C$ of $F$, both maps are just the standard isomorphism
$I\udot\cong\O$, so $\phi$ is invertible.
\end{proof}

\subsection{Obstruction theory}
In Section \ref{dcat}, the association of the complex 
$$I\udot= \big\{ \O_X \stackrel{s}{\rightarrow} F \big\} \in D^b(X)$$
to the pair $(F,s)$ was seen to be injective on objects.
By Theorem \ref{def}, the fixed determinant deformation theory in $D^b(X)$ matches the deformation theory of pairs to all orders.
Hence, $P(X)$ is a component{\footnote{We drop
the subscripts $n$ and $\beta\neq 0$ for notational convenience.
Formally, a component is used here to signify a union of connected
components.}}
 of the
moduli space of complexes of trivial determinant in $D^b(X)$. We will use the
obstruction theory  $\Ext^*(I\udot,I\udot)_0$ to 
define a
virtual fundamental class
on $P(X)$.

Since stable pairs have no nontrivial automorphisms, the moduli
space $P(X)$ is fine. 
There is a \emph{universal} stable pair{\footnote{
Pairs are better than 
stable
sheaves which have scalar automorphisms. In general, there is only
a universal \emph{twisted} sheaf on the product of $X$ and the
moduli space of stable sheaves.}} 
\beq{uni}
\O_{X\times P(X)}\to\FF 
\eeq
on $X\times P(X)$.
The moduli space is the GIT quotient of a subset of the product of 
a Quot scheme and a Grassmannian \cite{LPPairs1}.
There is a universal sheaf pulled back from the Quot scheme, with a universal
section over the product. Over the stable locus there are no stabilizers,
so by Kempf's lemma (see for instance 
\cite[Theorem 4.2.15]{HLShaves}) the universal pair descends to 
$P(X)$.

The universal pair \eqref{uni} determines a universal complex in
the derived category,
$$
\mathbb I\udot=\{\O_{X\times P(X)}\to\FF\} \in D^b(X\times 
P(X)),
$$
with $\FF$ flat over $P(X)$.
Let $\pi$ denote the projection 
$$\pi\colon X\times P(X) \to P(X).$$
We will prove that 
the complex 
$R\pi_*R\hom(\mathbb I\udot,\mathbb I\udot)_0$ determines
an obstruction theory on $P(X)$.

Since $X\times P(X)$ is projective and $\mathbb I\udot$ is perfect we may resolve it by a finite complex of locally 
free sheaves $A\udot$, and form
$$
(A\udot)^\vee\otimes A\udot\ \cong\ \O_{X\times P}\,\oplus\,
((A\udot)^\vee\otimes A\udot)_0.
$$
The first summand is the image of the identity map and the second is
the kernel of the trace map. These split each other since 
$$\tr\circ\id=\rk(A\udot)
=\rk(\mathbb I\udot)=1=\id\circ\tr.$$ 
We define the quasi-isomorphism class
of the trace-free Homs by
$$
R\hom(\mathbb I\udot,\mathbb I\udot)_0\simeq((A\udot)^\vee\otimes A\udot)_0.
$$
The following result 
is a standard consequence of Lemma \ref{homs} and the
Nakayama Lemma, but we give the argument in full.

\begin{lem} \label{perf} 
The complex $R\pi_*R\hom(\mathbb I\udot,\mathbb I\udot)_0$
on $P(X)$ is quasi-iso\-morphic
to a 2-term complex of locally free sheaves $\{E_1\to E_2\}$.
\end{lem}

\begin{proof}
Let $B\udot$ be
a sufficiently negative locally free resolution of the complex
$((A\udot)^\vee\otimes A\udot)_0$ trimmed to start at least
4 places earlier than $((A\udot)^\vee\otimes A\udot)_0$.
Then, by standard arguments, for all $j$,
\begin{enumerate}
\item[(i)]
$R^{\leq 2}\pi_* B^j = 0$,
\item[(ii)]
$R^3\pi_* B^j$ is locally free.
\end{enumerate}
The complex $E\udot$ with
$$E^k \cong R^3\pi_* B^{k+3}$$
is finite, locally free, and quasi-isomorphic to
$R\pi_*R\hom(\mathbb I\udot,\mathbb I\udot)_0$.


By base change,
the restriction of $E\udot$ over a point $[I\udot]\in
P(X)$ is a complex of vector spaces 
computing $\Ext^*(I\udot,
I\udot)_0$.     
By Lemmas \ref{homs} and \ref{dddd}, the local-to-global
spectral sequence, and Serre duality, $\Ext^i(I\udot,I\udot)$ 
is nonzero only for $i$ between 0 and 3.

By  \eqref{ddd} and the local-to-global spectral
sequence, we have the
isomorphisms 
\begin{eqnarray*}
K_X& \Rt{\id}& \hom(I\udot,I\udot\otimes K_X), \\
H^0(K_X)&\Rt{\id}&\Hom(I\udot,I\udot\otimes K_X).
\end{eqnarray*}
Then, by Serre duality,
 $$\Ext^3(I\udot,I\udot)\Rt{\tr}H^3(\O_X)$$
is also an isomorphism and $\Ext^3(I\udot,I\udot)_0$ vanishes.
Similarly, the composition 
$$\C\Rt{\id}\Hom(I\udot,I\udot)\Rt{\tr}H^0(\O_X)$$
is multiplication by $\rk(I\udot)=1$. Hence $\Hom(I\udot,I\udot)_0$
also vanishes and
the trace-free extensions $\Ext^i(I\udot,I\udot)_0$
are concentrated entirely in degrees $i=1$ and $2$. After base change
to any
point $[I\udot]\in P(X)$, the complex $E\udot$ has cohomology only
in degrees 1 and 2.

If $E^{n>2}$ is the last nonzero term of $E\udot$, 
then on each fiber 
\begin{equation}\label{xdee}
E^{n-1}\to E^n
\end{equation}
 is surjective. Hence, the map
 \eqref{xdee} is
surjective globally
with locally free kernel. Replacing $E^{n-1}$ by the kernel
and $E^n$ by zero, we can inductively assume $n=2$.
Similarly if the first nonzero term is $E^{m<1}$, 
then $E^m\to E^{m+1}$ is injective on fibers with
locally free cokernel.  We conclude $E\udot$ is quasi-isomorphic to a 2-term
complex $\{E_1\to E_2\}$.
\end{proof}

The Atiyah class of the universal complex $\mathbb I\udot$ gives an element
of
\begin{equation}\label{vgrt}
\Ext^1(\mathbb I\udot,\mathbb I\udot\otimes L\udot_{X\times P(X)}),
\end{equation}
where $L\udot_{X\times P(X)}$ is 
(a locally free resolution of) the
cotangent complex of $X\times P(X)$. 
Since the cotangent complex of
a product is the sum of the (pullbacks of the) cotangent complexes
of the factors, the
$\Ext$ group \eqref{vgrt}
maps to
$$
\Ext^1(R\hom(\mathbb I\udot,\mathbb I\udot)_0,\pi^*L
\udot_{P(X)}),
$$
which by Serre duality along $\pi$ is isomorphic to
$$
\Ext^{-2}(R\pi_*(R\hom(\mathbb I\udot,\mathbb I\udot)_0\otimes\omega_\pi),
L\udot_{P(X)}).
$$
The relative dualizing sheaf $\omega_\pi$ is the pull-back of the canonical
bundle of $X$. We obtain a map
\beq{obsthy}
R\pi_*(R\hom(\mathbb I\udot,\mathbb I\udot)_0\otimes\omega_\pi)[2]\to
L\udot_{P(X)}.
\eeq

We claim that \eqref{obsthy} 
provides a perfect 
obstruction theory for $P(X)$ in the sense of \cite{BFNormalCone} 
via the results of \cite{HT,LieblichModuli,Lowen}. 
More precisely,
the obstruction classes of
\cite{LieblichModuli,Lowen} are expressed in \cite{HT} as a product of Atiyah
and Kodaira-Spencer 
classes. As a result, compatibility
with \eqref{obsthy} is 
obtained.{\footnote{See Theorem 4.1 \cite{HT} for a careful discussion of
the obstruction class and the results needed to
prove that \eqref{obsthy} is an obstruction theory for $P(X)$.}}
By Lemma \ref{perf} and Serre duality, the obstruction theory \eqref{obsthy}
is \emph{perfect}: representable by a 2-term complex of locally free sheaves
$$E_2^\vee\to E_1^\vee$$ in degrees $-1$ and 0
respectively. Therefore, the results of 
\cite{BFNormalCone, LiTianVirtual} yield a virtual class.

\begin{thm} 
$P_{n}(X,\beta)$ carries an algebraic virtual class
$$
[P_{n}(X,\beta)]^{vir} \in A_{c\_\beta}(P_n(X,\beta),\Z)
$$
where
$$
c\_\beta=\int_\beta c_1(X)=-\chi(R\Hom(I\udot,I\udot)_0)
$$
is the virtual dimension of $P_{n}(X,\beta)$ 
with the obstruction theory inherited from the moduli space of fixed determinant
complexes $I\udot\in D^b(X)$.
\end{thm}

Consider a smooth family of projective 3-folds,
$$\pi:\mathcal X\to B,$$
over a base $B$ with central fiber
$X \cong \mathcal
X\times_B\{0\}$.
 Let $$\mathcal P_n(\mathcal X,\beta)\to B$$ denote
the relative moduli space of stable pairs on the fibers of $\pi$, 
and let
$$
i_0\colon P_n(X,\beta)\to\mathcal P_n(\mathcal X,\beta)
$$
be the inclusion of the space of stable pairs on the central fiber $X$. 
The 
deformation invariance of the virtual class in the
following form is obtained
from
\cite[Corollary 4.3]{HT}.

\begin{thm} 
There is $\pi$-relative virtual class
$$
[\mathcal P_{n}(\mathcal X,\beta)]^{vir}\in A_{c\_\beta+\dim B}(\mathcal
P_n(\mathcal X,\beta),\Z)
$$
for which
$$
i_0^![\mathcal P_{n}(\mathcal X,\beta)]^{vir}=[P_{n}(X,\beta)]^{vir}.
$$
\end{thm}

\subsection{Invariants}
If $X$ is a  Calabi-Yau 3-fold, $c\_\beta=0$ for every curve
class $\beta \in H_2(X,\Z)$.

\begin{defn} If $c\_\beta=0$, the stable pairs invariants 
$P_{n,\beta}\in\Z$ are defined
to be the degree of the virtual cycle 
$$P_{n,\beta} = \int_{[P_n(X,\beta)]^{vir}} 1.$$
\end{defn}

We define similar invariants with primary field
insertions in 
case $c\_\beta>0$ in Section \ref{fano}.

For future reference, we list the letters we use for the various curve counting
invariants associated to the class $\beta\in H^2(X,\Z)$. 
The definitions and references will be given later.
$$\begin{array}{lll}
\bullet \  N_{g,\beta}& : &\text{genus $g$ connected Gromov-Witten 
invariant.}\\
\bullet \  N^\bullet_{g,\beta}&: &\text{genus $g$ disconnected 
Gromov-Witten invariant} \\ & & \text{with no contracted contributions.} \\
\bullet \  n_{g,\beta}&: & \text{genus $g$ Gopakumar-Vafa BPS invariant.} \\
 \bullet \ I_{n,\beta}&: & \text{DT invariant with 
Euler characteristic $n$.}\\
\bullet\  P_{n,\beta}&:& \text{stable pairs invariant with 
Euler characteristic $n$.}
\end{array}$$

\section{Conjectures} \label{conjs}

\subsection{GW/DT for Calabi-Yau 3-folds}
We recall the conjectural GW/DT correspondence
from \cite{MNOP1} for a Calabi-Yau 3-fold $X$.

The disconnected Gromov-Witten invariants of $X$
for nonzero curve classes $\beta\in H_2(X,\Z)$ are
$$N^\bullet_{g,\beta} = \int_{[\overline{M}^\bullet_{g}(X,\beta)]^{vir}} 1,$$
where $\overline{M}^\bullet_{g}(X,\beta)$ is the moduli space
of stable maps with
possibly disconnected domains and  {\em no} contracted connected
components.{\footnote{We
follow here the notation of \cite{BryanP}.}}
Let 
$$
Z_{GW,\beta}(u)=\sum_g N_{g,\beta}^\bullet \ u^{2g-2}.
$$
be the partition function.
Alternatively, the partition function may be defined via the exponential 
of the connected potential,
$$
Z_{GW}(u,v)=1+\sum_{\beta\ne0}Z_{GW,\beta}(u) v^\beta=\exp F_{GW}(u,v),
$$
where
$$
F_{GW}(u,v)=\sum_{\beta\ne0}
\sum_g N_{g,\beta}\ u^{2g-2}v^\beta,
$$
and
$$N_{g,\beta} = \int_{[\overline{M}_{g}(X,\beta)]^{vir}} 1$$
is the standard connected Gromov-Witten invariant.

Let $I_n(X,\beta)$ be the Hilbert scheme of subschemes 
$$Z\subset X$$
with holomorphic Euler characteristic $n$ and fundamental class 
$\beta$.
To obtain a virtual class \cite{ThCasson}, the Hilbert scheme is viewed as
a moduli space of ideal sheaves.{\footnote{As in the case of stable
pairs, the straightforward obstruction theory of the Hilbert
scheme is not appropriate and an alternative is required.}} 
The DT invariants are defined by
$$I_{n,\beta} = \int_{[I_{n}(X,\beta)]^{vir}} 1.$$
The partition function of DT
theory is 
$$
Z_{DT,\beta}(q)=\sum_n I_{n,\beta}\ q^n.
$$
By the boundedness of the Hilbert scheme $I_n(X,\beta)$, the number of free and embedded points of any element $Z$ ---
the length of the maximal 0-dimensional subsheaf of $\O_Z$ --- 
is bounded above.
Therefore $I_n(X,\beta)$ is empty for $n$ sufficiently negative, and $Z_{DT,\beta}(q)$ is a Laurent series in $q$.

The irreducible components of the subscheme $Z$ consist of curves and 
0-dimensional
subschemes which wander all over $X$.  
A \emph{reduced}
partition function is defined \cite{MNOP1}
by dividing out by the degree 0 series,
$$
Z'_{DT,\beta}(q)=\frac{Z_{DT,\beta}(q)}{Z_{DT,0}(q)}\,.
$$
The degree 0 series is evaluated by the formula
$$Z_{DT,0}(q)=M(-q)^{\chi(X)},$$
conjectured in \cite{MNOP1} and proved in \cite{BFHilb,LPCobordism,LiDT}.
Here,
$$M(q)=
\prod_{n\ge1}(1-q^n)^{-n}$$ is the MacMahon function.

The reduced series
$Z'_{DT,\beta}(q)$ is conjectured in \cite{MNOP1} to be the 
Laurent expansion
of a rational function in $q$ 
invariant under the transformation $q\leftrightarrow q^{-1}$.
The GW/DT correspondence of \cite{MNOP1} is the conjectural
equality
$$
Z'_{DT,\beta}(q)\ =\ Z_{GW,\beta}(u)
$$
after the variable change $-q=e^{iu}$.

\subsection{Stable pairs conjectures for Calabi-Yau 3-folds}
\label{vvvv}
Our new invariants $P_{n,\beta}$
counting stable pairs on $X$ do not have the drawback of freely
roaming points, 
so a
reduced partition function is not necessary.
For nonzero $\beta \in H_2(X,\Z)$, let
\beq{pairsgen}
Z_{P,\beta}(q)=\sum_n P_{n,\beta}\ q^n
\eeq
be the partition function of the stable pairs theory.
The moduli spaces $P_n(X,\beta)$ are
empty for the same sufficiently negative $n$ as for the $I_n(X,\beta)$, so
 $Z_{P,\beta}(q)$ is a Laurent series in $q$. 

\begin{conj}
\label{111} 
The partition function
$Z_{P,\beta}(q)$ is the 
Laurent expansion of a rational function in $q$
 invariant under $q\leftrightarrow q^{-1}$.
\end{conj}

In fact, the above rationality Conjecture will be significantly
refined after our discussion of the BPS state counts of
Gopakumar and Vafa.

\begin{conj} \label{222}
All the partition functions coincide,
$$Z_{P,\beta}(q)\ =\ Z'_{DT,\beta}(q)\ =\ Z_{GW,\beta}(u),$$
after the variable change $-q=e^{iu}$.
\end{conj}

It appears all reasonable enumerative theories of curves
on Calabi-Yau 3-folds
are actually equivalent.

\subsection{Wall crossing formula}
The first equality of Conjecture \ref{222},
\beq{conjj}
Z_{P,\beta}(q)\cdot Z_{DT,0}(q)\ = \
Z_{DT,\beta}(q)
\eeq
can be expanded to yield
\beq{wallnos}
\sum_m   P_{n-m,\beta}\cdot    I_{m,0}  =I_{n,\beta}.
\eeq
Relation \eqref{wallnos}
should be interpreted as a wall-crossing formula for
counting invariants in the derived category of coherent sheaves $D^b(X)$
under a change of stability condition \cite{BrStability}. So far, however,
counting invariants have yet to be defined for general complexes,
and Bridgeland stability conditions have not been shown to exist for compact 3-folds.

Let us assume that there is a Bridgeland stability
condition for which ideal sheaves $\I_Z$ of subschemes $Z\subset X$ 
satisfying
$$\chi(\O_Z)=n , \ \ [Z]=\beta$$
are stable and constitute the moduli space 
of semistable objects of the same phase and Chern character
$$
(\mathrm{rk},ch_1,ch_2,ch_3)=(1,0,-\beta,-n)
$$
and trivial determinant. 
The counting invariant here exists \cite{ThCasson} and is $I_{n,\beta}\in\Z$.
We also assume the structure sheaves of single points are stable. Since
structure sheaves of subschemes of length $m>1$ are only semistable (and
have automorphisms), it is not clear what their counting invariant should
be, but the virtual number of \emph{ideal} sheaves of such subschemes is
$I_{m,0}$.

Now move the central charge of \cite{BrStability} across
 a codimension
1 wall along which the phase of $\I_Z$ equals
 the phase of the structure sheaf of a point $\O_p$ minus $1$:
$$
\phi(\I_Z)=\phi(\O_p)-1=\phi(\O_p[-1]).
$$
Any free or embedded points 
$p$ in $Z$ give rise to exact sequences of the form
$$
0\to\I_Z\to\I_{Z'}\to\O_p\to0
$$
where $Z'$ is $Z$ with the point $p$ removed.
In $D^b(X)$, we obtain an exact triangle
$$
\O_p[-1]\to\I_Z\to\I_{Z'},
$$
destabilizing $\I_Z$ as we cross the wall. 

More generally, 
let $Q\subset\O_Z$ denote the maximal subsheaf with 0-dimensional 
support (roughly the
structure sheaf of the union of all free and embedded points). Then
\beq{destab}
Q[-1]\to\I_Z\to\I_{Z'}
\eeq
is the maximal destabilizing extension making up $\I_Z$. 
Conversely, previously unstable objects become stable as we cross the wall. 
These are objects $I\udot\in
D^b(X)$ which are extensions of 
the form (\ref{destab}) but in the opposite direction,
\begin{equation}\label{ccf}
\I_{Z'}\to I\udot\to Q[-1],
\end{equation}
classified by elements of $\Ext^2(Q,\I_{Z'})$. We recognize \eqref{ccf}
 as the
form of (\ref{b}). We expect that the
moduli space of pairs $P_{n}(X,\beta)$ 
gives precisely the space of stable objects for the new 
stability condition.{\footnote
{Recently, Bayer \cite{BayerPoly} and Toda \cite{TodaStab}
have defined variants of Bridgeland's axioms for a 
stability condition within which the above wall crossing
occurs.}}
The virtual numbers $P_{n,\beta}$ should then be
the right counting invariants.

Relation (\ref{wallnos}) 
has the form of a wall crossing formula 
envisaged by Joyce 
\cite{JoWallCrossing} for invariants counting stable objects in
$D^b(X)$. The formula expresses all
of the possible ways that a stable object on one side of the wall -- an ideal sheaf $\I_Z$ -- can be written as extensions
of objects in $D^b(X)$ which are stable on the other side of the wall. 
No worse configurations occur since the K-theory classes of $\I_Z$, $\I_{Z'}$
and $\O_p$ are all primitive and distinct.
So the $m^{th}$ term in (\ref{wallnos}) is the contribution from subschemes
$Z$ whose maximal 0-dimensional subscheme (or total number of free
and embedded points) is of length $m$ -- the length of $Q$. 
We expect that in Joyce's theory the space of extensions between the semistable pieces
$\I_{Z'}$ and $Q[-1]$ contributes
\begin{equation}\label{bbgg}
-\chi(\I_{Z'},Q[-1])=m.
\end{equation}
As $Z'$ and $Q$ vary this gets multiplied by the product of the invariants associated to $\I_{Z'}$ and $Q[-1]$. The former should be $P_{n-m,\beta}$. Formula
\eqref{wallnos} requires the latter to contribute $I_{m,0}/m$. That is,
the automorphisms of the structure sheaves of length-$m$ points 
should affect the virtual
count for ideal sheaves of points by the factor $1/m$.

Until Joyce's theory can be fully extended to the derived category and made
to include the virtual class, our entire discussion here
remains conjectural.

Similar phenomena have been observed recently in a noncommutative example
\cite{SzendroiDT}. The paper \cite{DenefMoore} also studies 
counting invariants (and more generally BPS states) and 
wall crossing in $D^b(X)$ from a physical point of view. 
In fact, Denef and Moore
 predict a wall crossing of precisely the form (\ref{conjj}) 
as the K\"ahler form and, crucially, the B-field, cross a certain wall in the
space of stability
conditions
(formulae (6.21)--(6.24) of \cite{DenefMoore}).
 Their predictions are based on a supergravity analysis rather
than algebraic geometry. For a certain stability condition,
the virtual number of
stable objects in $D^b(X)$ is predicted without
identifying the form of these objects.
Other wall crossings studied in \cite{DenefMoore,DiacMoore} 
have not yet been considered mathematically.

\subsection{BPS refinement}
Motivated by the M-theory prediction of Gopa\-kumar and
Vafa \cite{GV1, GV2} 
for Calabi-Yau 3-folds $X$, we conjecture
a much stronger rationality statement for $Z_{P,\beta}(q)$.

Let $F_P(q,v)$ denote the {\em connected} series
for the stable pairs invariants,
$$
F_P(q,v)
=\log Z_P(q,v)=\log\Big(1+\sum_{\beta\ne0}Z_{P,\beta}(q)v^\beta\Big).
$$
The functions $F_{P,\beta}$ are defined by
$$
F_P(q,v)=\sum_{\beta\ne0}F_{P,\beta}(q)v^\beta,
$$
While $F_{P,\beta}(q)$ should be thought of as counting pairs whose
support is connected,
we do not have a geometric method to define the invariants.
Though
$F_{P,\beta}(q)$ is still a Laurent series,
the coefficients of $F_{P,\beta}(q)$
need not be integral due to the logarithm.

Let $V_0\subset \mathbb{Q}(q)$ be the linear space 
of Laurent polynomials invariant
under
$q\leftrightarrow q^{-1}$. 
Certainly,
$$\Phi(q) = \frac{(1-q)^2}{q} = q -2 + \frac{1}{q}\in V_0.$$
\begin{lem}\label{ddumb}
Every Laurent series 
$L(q)\in \mathbb{Q}(\!(q)\!)$ can be {\em uniquely} written
as a sum
\begin{equation}\label{frre}
L(q) = \sum_{g> -\infty} l_g \ \Phi(-q)^{g-1}, \ \ l_g \in \mathbb{Q}
\end{equation}
with only {\em finitely} many
positive $g$ terms.
\end{lem}
\begin{proof}
Any $g>1$ term
$$
\Phi(-q)^{g-1}=\left((-q)-2+(-q)^{-1}\right)^{g-1} \in V_0$$
has a pole of order $g-1$ at $q=0$. 
Together the $g>1$ terms in \eqref{frre}
specify the finite polar part of $L(q)$.
The $g\leq 1$ term
$$\Phi((-q))^{g-1} = \left(\frac{(-q)}{(1-(-q))^2}\right)^{1-g} 
= (-q)^{1-g} + 2(1-g)(-q)^{2-g} + \ldots$$
is regular at $q=0$ with a zero of order $1-g$.
The right side of \eqref{frre} is thus uniquely determined.
\end{proof}

For nonzero $\beta\in H_2(X,\Z)$, let $\div(\beta)$
denote the divisibility of $\beta$.
For simplicity, we either assume $H_2(X,\Z)$ is torsion-free or 
allow $\beta$ take values in the torsion-free quotient
$$H_2(X,\Z)/\tau(H_2(X,\Z))$$ so each
stable pairs invariant is a sum of all classes differing 
from $\beta$ by torsion.{\footnote{In fact, no such assumption is necessary.
The full torsion information can be kept. 
Then, the interior sum in Lemma \ref{BPS} is over all
elements of $H_2(X,\Z)$ whose $r^{th}$ multiple is
$\beta$.}}

\begin{lem} The set of equations 
\label{BPS}
\begin{equation*}
\left\{ F_{P,\beta}(q)\ 
=\ \sum_{g> -\infty}\ \sum_{r|\div(\beta)}
n_{g,\frac{\beta}r}\frac{(-1)^{g-1}}r
\left((-q)^r-2+(-q)^{-r}\right)^{g-1} \right\}_{\beta\neq 0},
\end{equation*}
has a unique solution 
$\{n_{g,\beta}\}_{g>-\infty,\,\beta\neq 0}$. 
\end{lem}

\begin{proof}
We proceed inductively on the divisibility of $\beta$. If
$\div(\beta)=1$, then the statement is 
a direct consequence of Lemma \ref{ddumb} applied to $F_{P,\beta}(q)$.
For the inductive step, we apply Lemma \ref{ddumb}
to
$$F_{P,\beta}(q) -\ \sum_{g> -\infty}\ \sum_{r|\div(\beta),\ r>1}
n_{g,\frac{\beta}r}\frac{(-1)^{g-1}}r
\left((-q)^r-2+(-q)^{-r}\right)^{g-1}$$
to conclude the existence and uniqueness
of the solution.
\end{proof}

From vanishing of Lemma \ref{ddumb} applied inductively, we obtain
the following result.

\begin{lem}
For fixed $\beta$, 
$n_{g,\beta} = 0 $
for all sufficiently large $g$.
\end{lem}

We {\em define}
 the Gopakumar-Vafa
BPS state counts in genus $g$ and class $\beta\neq 0$ 
for $X$ to be the solutions $n_{g,\beta}$ of Lemma \ref{BPS}.
Our motivation for the
definition
is the agreement of the formula
$$F_{P}(q,v) = \sum_{g > -\infty}\sum_{\gamma\neq 0} 
\sum_{r\geq 1} n_{g,\gamma}\frac{(-1)^{g-1}}{r}
\ \Phi((-q)^r)^{g-1} v^{r\gamma},$$ 
after truncation by Conjecture \ref{333} below and the
variable change $$-q=e^{iu},$$ with the string theoretic
Gopakumar-Vafa formula 
\begin{equation}\label{bttr}
F_{GW}(u,v) = \sum_{g\geq 0} \sum_{\gamma\neq 0} n_{g,\gamma}\ u^{2g-2}
 \sum_{r\geq 1} \frac{1}{r} 
\left(\frac{\sin(ru/2)}{u/2}\right)^{2g-2} v^{r\gamma}
\end{equation}
for Gromov-Witten theory via BPS counts.

Philosophically, 
we view the integrals $P_{n,\beta}$ as giving a rigorous
treatment of the heuristic approach{\footnote{S. Katz \cite{KatzBPS}
has proposed a mathematical  definition of the $g=0$ BPS state counts
via virtual Euler characteristics of moduli of
semistable sheaves on curves in $X$. Katz's construction, in which
neither pairs nor derived categories appear, is nevertheless
not unrelated to our perspective.}}
 to BPS
state counting proposed in \cite{KKVSpinning} 
via virtual Euler characteristics of
 Hilbert schemes of points on curves in $X$.
Certainly, $P_n(X,\beta)$ may be thought of as
a compactification of the Hilbert scheme of points
on curves in $X$. And,
$P_{n,\beta}$ is precisely{\footnote{In the Calabi-Yau
case, the obstruction theory on $P_n(X,\beta)$ is
self-dual \cite{BehrendDT,DT} 
and the virtual count $P_{n,\beta}$ is the
virtual 
Euler characteristic, up to sign.}} a signed virtual Euler characteristic.

We replace the rationality statement of Conjecture \ref{111} with
a much stronger and more geometric vanishing.

\begin{conj} The invariants $n_{g,\beta}$ vanish for $g<0$.
\label{333}
\end{conj}

Define the linear subspace $V_d\subset \mathbb{Q}(q)$ for $d>0$
by the following spanning set:
\beq{Vd}
V_d = \text{Span}_{\mathbb{Q}}\big\{ \Phi((-q)^r)^{g-1} 
\big\}_{g\geq 0, \ 1\leq r \leq d}\,.
\eeq
Since $V_0$ is the linear space 
of Laurent polynomials invariant
under $q\leftrightarrow q^{-1}$, $V_0$ is
 spanned by $\{\Phi(-q)^{g-1}\}_{g\ge1}$
and  contained in all $V_d$.
A basis of $V_d/V_0$ is given by the $g=0$ terms of \eqref{Vd}, 
$$\frac{(-q)}{(1-(-q))^2}\,,\ \frac{(-q)^2}{(1-(-q)^2)^{2}}\,,\ \ldots,\
\frac{(-q)^d}{(1-(-q)^d)^2}\,.$$

Conjecture \ref{333}  implies the nontrivial inclusion
\begin{equation}\label{tt55}
F_{P,\beta}(q) \in V_{\div(\beta)}.
\end{equation}
As a consequence, we conclude that
both $F_{P,\beta}(q)$ and $Z_{P,\beta}(q)$ 
have possible poles {\em only at  $r^{th}$ roots of unity}
where $r\leq\div(\beta)$. 

In fact, inclusion \eqref{tt55} is much
stronger than
the statement about poles. For example,
$$\Psi(q) = \frac{q}{1+q+q^2}$$
is invariant under $q\leftrightarrow q^{-1}$, has poles at only
$3^{rd}$ roots of unity,
but is easily seen to satisfy
$$\forall d, \ \ \Psi(q) \notin V_d.$$

Conjecture \ref{333} implies the following two effectivity
statements. The asterisk denotes the dependence on
the conjecture.

\begin{lems} \label{34r}
$F_{P,\beta}(q)$ is uniquely and
effectively determined
by the coefficients of order $q^n$ for $-\infty<n\leq\div(\beta)$.
\end{lems}

\begin{proof}
We can uniquely write $F_{P,\beta}$ as an element of $\Theta(q)\in V_0$
plus a linear combination of
\begin{equation}\label{nnnww}
\frac{(-q)}{(1-(-q))^2}\,,\ \frac{(-q)^2}{(1-(-q)^2)^{2}}\,,\ \ldots,\ 
\frac{(-q)^{\div(\beta)}}{(1-(-q)^{\div(\beta)})^2}\,.
\end{equation}
The functions \eqref{nnnww} are regular at $q=0$.
Therefore,
the coefficients of $q^n$ for $n\leq 0$ uniquely determine $\Theta(q)$,
and the coefficients of $q^n$ for $1\leq n \leq \div(\beta)$
determine the linear combination \eqref{nnnww}.
\end{proof}

\begin{lems} \label{gg12}
$Z_{P}(q,v)$ is uniquely and
effectively determined
by the invariants 
$
\big\{ P_{n,\beta}\colon-\infty<n\leq 1,\,\forall\beta\big\}$.
\end{lems}

\begin{proof}
The argument follows the proof of Lemma \ref{34r}.
Since data is given for all classes $\beta$ simultaneously, we may work by
induction on the degree of $\beta$. The coefficients of $v^\beta$ when
$\beta$ is primitive are dealt with by Lemma \ref{34r}. For general $\beta$
we need only handle the $r=1$ terms in 
\begin{equation*}
 F_{P,\beta}(q)\ 
=\ \sum_{g\ge0}\ \sum_{r|\div(\beta)}
n_{g,\frac{\beta}r}\frac{(-1)^{g-1}}r
\left((-q)^r-2+(-q)^{-r}\right)^{g-1},
\end{equation*}
since the others are dealt with by the induction assumption.
\end{proof}

The proof of Lemma$^*$ \ref{gg12} shows
 the statement
is local in the following sense. To determine $Z_{P,\gamma}(q)$,
the data $
\big\{ P_{n,\beta}\big\}_{n\leq 1}$
is required only for classes
$\beta$ which are effective summands of $\gamma$.

\subsection{BPS integrality}

We show the integrality
constraints on the stable pairs counts $P_{n,\beta}$ exactly
matches the BPS integrality of $n_{g,\beta}$.
We do not require Conjecture \ref{333} for the result.

\begin{thm} The invariants
 $P_{n,\beta}$ are integers for all $n,\beta\neq 0$
 if and only if
the coefficients $n_{g,\beta}$ are integers for all $g,\beta\neq 0$.
\end{thm}

\begin{proof}
By definition,
$Z_{P,\beta}$ is the $v^\beta$ coefficient of
$$
\exp\left( \sum_{g>-\infty} \sum_{\gamma\neq 0}
\sum_{k\geq 1} n_{g,\gamma}\frac{(-1)^{g-1}}k
\left((-q)^k-2+(-q)^{-k}\right)^{g-1}v^{k\gamma}\right).
$$
To simplify notation, let
$$\tilde n_{g,\gamma}=(-1)^{g-1}n_{g,\gamma},\ \ Q=-q.$$
Then, $Z_{P,\beta}$ is the $v^\beta$ coefficient of
\beq{expanded}
\exp\left( \sum_{g>-\infty} \sum_{\gamma\neq 0}
\sum_{k\geq 1} \frac{\tilde{n}_{g,\gamma}}{k}
\Phi(Q^k)^{g-1} v^{k\gamma}\right).
\eeq

We have seen already in the proof of Lemma \ref{ddumb} that for all $k\in \Z$,
\beq{hots}
\Phi(Q)^k=\sum_{l\geq -k} \phi^{k}_l Q^l, \ \ \phi^{k}_l \in \mathbb{Z},
\eeq
with leading coefficient $\phi^{k}_{-k}=1$.

Let $L$ be a very ample line bundle on $X$.
We prove the Theorem inductively on the degree 
$$L_\beta = \int_\beta c_1(L)$$ of 
$\beta$. 

Classes $\beta \neq 0$ of minimal degree are necessarily primitive.
Hence, only terms 
with $\gamma=\beta$ and $k=1$ 
contribute to the  $v^\beta$ coefficient of \eqref{expanded},
$$
Z_{P, \beta}=\sum_{g>-\infty} \tilde n_{g,\beta}\Phi(Q)^{g-1}.
$$
By the expansion \eqref{hots},
the coefficients of $Z_{P,\beta}$ and the invariants
$\tilde n_{g,\beta}$
are related by a triangular transformation
 with 1s along the diagonal.
Hence,
$$\text{$\big\{ P_{n,\beta}\big\}_{n\in \Z}$ are integral $\iff$
$\big\{\tilde n_{g,\beta}\big\}_{g\in \Z}$ are integral.}$$

For the induction step, we write $Z_{P,\beta}(q)$ with the
$\gamma=\beta$ and $k=1$ term in front:
\begin{equation}\label{mess}
\sum_{g>-\infty} \tilde n_{g,\beta}\Phi(Q)^{g-1}\ 
+ \left(\exp\
\sum_{g>-\infty}\sum_{\gamma\ne\beta} \sum_{k\geq 1}
\frac{\tilde n_{g,\gamma}}k \Phi(Q^k)^{g-1}v^{k\gamma}\right)
\!\!\!\!\!\begin{array}{c} \\\\ v^\beta\end{array},
\end{equation}
where the suffix denotes taking the $v^\beta$ coefficient.
The term inside the brackets in
(\ref{mess}) is
$$
\exp\ \sum_{g>-\infty}\sum_{\gamma\ne\beta}\sum_{k\geq 1}
\sum_{l \ge1-g}
\frac{\tilde n_{g,\gamma}}k \phi^{g-1}_{l}Q^{kl}v^{k\gamma}.
$$
Doing the sum over $k$ first, we obtain
$$
\exp\ \sum_{g>-\infty}\sum_{\gamma\ne\beta}
\sum_{l \ge1-g}
{\tilde n_{g,\gamma}}\ \phi^{g-1}_{l}\left(-\log(1-Q^{l}v^{\gamma})\right).
$$
Substitution in \eqref{mess} yields the following expression
for $Z_{P,\beta}(q)$:
\begin{equation}\label{mess2}
\sum_{g>-\infty} \tilde n_{g,\beta}\Phi(Q)^{g-1}\ 
+
\left(
\prod_{g>-\infty}\prod_{\gamma\ne\beta} \prod_{l \ge1-g}
\!\!\!\left(\frac1{1-Q^lv^\gamma}\right)^{\!\!\tilde n_{g,\gamma}\cdot 
\phi^{g-1}_{l}}\right)
\!\!\!\!\!\begin{array}{c} \\\\ v^\beta\end{array}.
\end{equation}
By the argument in the base case, if all
the $\tilde n_{g,\gamma}$ contributing to the second term are integral,
 then 
$$\text{$\big\{ P_{n,\beta}\big\}_{n\in \Z}$ are integral $\iff$
$\big\{\tilde n_{g,\beta}\big\}_{g\in \Z}$ are integral.}$$
But all the $\tilde{n}_{g,\gamma}$ contributing to the second term correspond
to classes $\gamma$ of strictly lower degree than $\beta$,
so the induction hypothesis applies.
\end{proof}

Since the stable pair invariants $P_{n,\beta}$ are certainly integral,
we conclude the BPS counts defined by $n_{g,\beta}$ are also integral.
Conversely, the integrality of $n_{g,\beta}$ imposes no further
conditions on $P_{n,\beta}$.

If the vanishing of Conjecture \ref{333} is assumed, we conclude
the integrality placed on Gromov-Witten theory by Conjecture \ref{222}
exactly coincides with the integrality predicted by the 
Gopakumar-Vafa formula \eqref{bttr}.
\label{vvvvv}
\subsection{The Fano case} \label{fano}
For an arbitrary  $3$-fold $X$, the virtual dimension of the moduli
space of pairs is
$$\text{dim}_{\mathbb{C}} [ P_n(X,\beta)] ^{vir} = c\_\beta,$$
where
$$c\_\beta=\int_\beta c_1(X).$$
For a nontrivial theory, we must have $c\_\beta \geq 0$.
In the {\em local Calabi-Yau} case satisfying
 $c\_\beta=0$, the discussion
is identical to the global Calabi-Yau case treated in Sections 
\ref{vvvv}-\ref{vvvvv}.
In the {\em local Fano} case satisfying  $c\_\beta>0$,
insertions in the theory are necessary.
We explain here the conjectural structure of the local Fano
case with primary field insertions.{\footnote{There exists
a full descendent theory which we will treat elsewhere.}}

Let $T_1,\ldots,T_m\in H^*(X,\Z)$ be a basis of the cohomology{\footnote{For
simplicity, we assume there is no odd cohomology. There
is no difficulty including the odd cohomology with
supercommuting variables.}} of $X$ mod torsion.
Let
$$N^\bullet_{g,\beta}(T_1^{e_1}\ldots T_m^{e_m})=
\int_{[\overline{M}^\bullet_{g,\sum_i e_i}(X,\beta)]^{vir}}
\prod_{i=1}^m \prod_{j=1}^{e_i} \text{ev}_{i,j}^*(T_i)$$
denote the 
disconnected Gromov-Witten invariant
with primary field insertions and no degree 0 connected components.
The primary fields should be viewed as simple incidence conditions.
The curves are required to intersect
the Poincar\'e duals of the insertions $T_i$.

For nonzero $\beta \in H_2(X,\Z)$, let
$$
Z_{GW,\beta}(u,t_i)=\sum_{g} \sum_{e\_\bullet}
N^\bullet_{g,\beta}(T_1^{e_1}\ldots T_m^{e_m})
\frac{t_1^{e_1}\ldots t_m^{e_m}}{e_1!\ldots e_m!} u^{2g-2}.
$$ 
The sum is over finitely many negative
and infinitely many positive $g$.
The second sum is over all vectors $(e_1,\ldots,e_m)$
of non-negative integers. 
Let
$$
Z_{GW}(u,v,t_i)=1+\sum_{\beta\ne 0}Z_{GW,\beta}(u,t_i)v^\beta
$$ 
be the full partition function.

Alternatively, $Z_{GW}=\exp F_{GW}$ is the exponential of the 
connected Gromov-Witten series,
$$
F_{GW}(u,v,t_i)=\sum_{\beta\ne0}F_{GW,\beta}(u,t_i)v^\beta,
$$ where
$$
F_{GW,\beta}(u,t_i)=\sum_{g\geq 0} \sum_{e\_\bullet}
 N_{g,\beta}(T_1^{e_1}\ldots T_m^{e_m})
\frac{t_1^{e_1}\ldots t_m^{e_m}}{e_1!\ldots e_m!}u^{2g-2}.
$$

To define the corresponding series via stable pairs, 
let 
$$\FF \rightarrow X\times P_{n}(X,\beta)$$
 denote the universal sheaf \eqref{uni}.
For a pair $$[\O_X\to F]\in P_{n}(X,\beta),$$ the restriction of
$\FF$
to the fiber
 $$X \times [\O_X \to F] \subset 
X\times P_{n}(X,\beta)
$$
is canonically isomorphic to $F$.
Let
$$\pi_X\colon X\times P_{n}(X,\beta)\to X,$$
$$\pi_P\colon X\times P_{n}(X,\beta)
\to P_{n}(X,\beta)$$
 be the projections on the first and second factors.
By definition, the operation
$$
\pi_{P*}\big(\pi_X^*(T_i)\cdot \text{ch}_2(\FF)
\cap(\pi_P^*(\ \cdot\ ))\big)\colon 
H_*(P_{n}(X,\beta))\to H_*(P_{n}(X,\beta))
$$
is the action of the primary field $\tau_0(T_i)$.

For nonzero $\beta\in H_2(X,\Z)$,
define the stable pairs invariant with primary field insertions by
\begin{eqnarray*}
P_{n,\beta}(T_1^{e_1}\ldots T_m^{e_m})&  = &
\int_{[P_{n}(X,\beta)]^{vir}}
\prod_{i=1}^m \tau_0(T_i)^{e_i} \\
& = & 
\int_{P_n(X,\beta)} \prod_{i=1}^m \tau_0(T_i)^{e_i}
\big( [P_{n}(X,\beta)]^{vir}\big).
\end{eqnarray*}
The partition function is 
$$
Z_{P}(q,v,t_i)=1+\sum_{\beta\ne0}Z_{P,\beta}(q,t_i)v^\beta,
$$
where
$$
Z_{P,\beta}(q,t_i)=\sum_{n} \sum_{e\_\bullet}
P_{n,\beta}(T_1^{e_1}\ldots T_m^{e_m})
\frac{t_1^{e_1}\ldots t_m^{e_m}}{e_1!\ldots e_m!}q^n.
$$
Since $P_n(X,\beta)$ is empty for sufficiently negative
$n$, $Z_{P,\beta}(q,t_i)$ is a Laurent series in $q$.
The connected stable pair invariants are obtained formally via
the logarithm,
$$
F_P(q,v,t_i)=\sum_{\beta\ne0}F_{P,\beta}(q,t_i)v^\beta
=\log Z_P(q,v,t_i),
$$
as in the Calabi-Yau case.

We state the rationality conjecture for $Z_{P,\beta}$ in
BPS form following the Gromov-Witten structure explained in
\cite{PandDegen,PandICM}.
The equation
\begin{multline}\label{nnd}
F_{P,\beta}(q,t_i) =
\sum_{g>-\infty} 
n_{g,\beta}(T_1^{e_1}\ldots T_m^{e_m})
\frac{t_1^{e_1}\ldots t_m^{e_m}}{e_1!\ldots e_m!} 
\\ \cdot {(-1)^{g-1}}\left((-q)-2+(-q)^{-1}\right)^{g-1}
\big(1+q\big)^{c\_\beta}.
\end{multline}
uniquely determines the invariants
$n_{g,\beta}(T_1^{e_1}\ldots T_m^{e_m})\in \mathbb{Z}$. 
The latter are
{\em defined} to be the BPS state counts. Moreover,
for fixed $\beta$,
$$n_{g,\beta}(T_1^{e_1}\ldots T_m^{e_m})=0$$
for all sufficiently large $g$.
The proofs are identical to those in the Calabi-Yau case.

\begin{conj}\label{444}
The invariants $n_{g,\beta}(T_1^{e_1}\ldots T_m^{e_m})$
vanish for $g<0$.
\end{conj}

Conjecture \ref{444} 
is the stronger form of
rationality obtained from the BPS perspective. As before,
we obtain an effectivity statement as a consequence.

\begin{lems}
$Z_{P}(q,v,t_i)$ is uniquely and
effectively determined
by the invariants 
$
\big\{ P_{n,\beta}(T_1^{e_1}\ldots T_m^{e_m})
\big\}_{n\leq 1}$.
\end{lems}

\begin{conj} \label{555} All the partition functions coincide
\begin{eqnarray*}
(-q)^{-\frac{c\_\beta}{2}} Z_{P,\beta}(q,v,t_i) &= & 
 (-q)^{-\frac{c\_\beta}{2}}
Z'_{DT,\beta}(q,v,t_i) \\
& = & (-iu)^{c\_\beta} Z_{GW,\beta}(u,v,t_i)
\end{eqnarray*}
after the change of variables $-q=e^{iu}$.
\end{conj}

The second equality in Conjecture \ref{555} is
the GW/DT correspondence for primary
fields in the local Fano case \cite{MNOP2}. We refer the reader to
\cite{MNOP2} for the definitions of the reduced partition function
$Z'_{DT,\beta}(q,v,t_i)$.

\subsection{Variants}
Let $G$ be a linearized algebraic group action on $X$.
The construction of the virtual class for the moduli
space of stable pairs in Section \ref{defsec} is valid
in the equivariant setting,
$$[P_n(X,\beta)]^{vir} \in A^G_{c\_\beta}(P_n(X,\beta),\mathbb{Z}).$$
We may then define an equivariant theory of stable pairs and
an equivariant correspondence following \cite{BryanP}. We leave
the details to the reader.

Another standard direction is the relative theory. While relative
Gromov-Witten theory has well-developed foundations 
 \cite{IonelParker,LiRuan,LiRelative, LiRelativeGW}, relative DT
theory awaits a definitive treatment. A sketch of the relative
theory for ideal sheaves (following suggestions of J. Li) is
given in \cite{MNOP2}. A very similar discussion holds for
the theory of stable pairs.

Let $X$ be a nonsingular projective 3-fold, and let 
$S\subset X$ be a nonsingular divisor.
The moduli space 
$P_{n}(X/S,\beta)$
parameterizes
stable relative pairs
\begin{equation}\label{vyq}
\O_{X[k]} \stackrel{s}{\rightarrow} F
\end{equation}
on $X[k]$, the $k$-step degeneration \cite{LiRelativeGW}
along $S$. $F$ is a sheaf on $X[k]$ with
$$\chi(F)=n$$
and whose support pushes down to the class
$$\beta\in H_2(X,\Z).$$
The stability conditions for the data are more complicated
in the relative geometry:
\begin{enumerate}
\item[(i)] $F$ is pure with finite locally free resolution,
\item[(ii)] the higher derived functors of the
restriction of $F$ -- to the singular loci of $X[k]$ and to the relative divisor $S_\infty\subset X[k]$ -- vanish,
\item[(iii)] the section $s$ has 0-dimensional cokernel supported
away from the singular loci of $X[k]$.
\item[(iv)] the pair \eqref{vyq} has only finitely many automorphisms covering
the automorphisms of $X[k]/X$.
\end{enumerate}

The moduli space $P_n(X/S,\beta)$ 
is a complete Deligne-Mumford stack
equipped with a map to the Hilbert scheme of points of $S$
via the relative geometry.
We expect
a perfect obstruction theory of virtual dimension 
$\int_\beta c_1(X)$ to be induced from the deformation theory of
complexes. 
All of these topics would benefit from
foundational work.

The relative theory of stable pairs should
admit a degeneration formula \cite{LiWuDegenDT} and a correspondence to relative
Gromov-Witten and relative DT theory following
\cite{MNOP2}.

Several other variants can be considered: families invariants for $3$-folds,
equivariant relative invariants, residue invariants in the presence of
a torus action, and so on. We expect the stable pairs theory
to be equivalent
in all reasonable cases to
the corresponding Gromov-Witten and DT
theories.

\section{First examples} \label{examples}
\subsection{Local $\PP^1$}
The simplest possible example is to consider the local Calabi-Yau 3-fold
$X$ given by the total space of the bundle 
$$\O_{\PP^1}(-1)^{\oplus2} \rightarrow \PP^1$$
in the class $[\PP^1]$ of the zero section.
Then, $P_{n}(X,[\PP^1])$
parameterizes nonzero sections, up to scale, of $\O_{\PP^1}(n-1)$ supported
on $\PP^1$. Thus,
$$
P_{n}(X,[\PP^1])\cong \text{Sym}^{n-1}(\PP^1)\cong\PP^{n-1}
$$
is nonsingular. As noted earlier, in the Calabi-Yau case, the
obstruction theory of the moduli space of stable pairs is self-dual.
Hence, if the moduli space is nonsingular,  
 the obstruction bundle is the cotangent bundle. So
$$
P_{n,[\PP^1]}=(-1)^{n-1}\chi_{top}(\PP^{n-1})=(-1)^{n-1}n,
$$
for $n\ge1$, and is 0 otherwise. Therefore
$$
Z_{P,[\PP^1]}(q)=q-2q^2+3q^3-\ldots=\frac{q}{(1+q)^2}\,,
$$
in complete agreement with 
$Z_{GW,[\PP^1]}$  and $Z'_{DT,[\PP^1]}$ \cite{FaberPand, MNOP1}.

For curve class $\beta=2[\PP^1]$, 
the lowest possible value of the holomorphic
Euler characteristic is 3. The moduli space $P_{3}(X,2[\PP^1])$ 
is a copy
of $\PP^1$ corresponding to the choice of a sub-bundle
$$
\O_{\PP^1}(-1)\subset\O_{\PP^1}(-1)^{\oplus2}.
$$
The associated stable pair is the structure sheaf of the curve
obtained by doubling along the sub-bundle with the canonical section.
So the potential starts as
$$Z_{P,2[\PP^1]}(q)=-2q^3+\ldots \, . $$

The next moduli space $P_{4}(X,2[\PP^1])$ is more interesting. Most stable
pairs correspond to a choice of sub-bundle
$$
\O_{\PP^1}(-2)\subset\O_{\PP^1}(-1)^{\oplus2}.
$$
We then take the structure sheaf of the doubling of the zero
section 
along the sub-bundle with the canonical section. 
We obtain an open set in
$\PP(H^0(\O_{\PP^1}(1))^{\oplus 2})\cong\PP^3$ consisting of pairs of sections
of $\O_{\PP^1}(1)$ which are not proportional to each other. 
Along the quadric
where the sections become proportional, we find instead an $\O_{\PP^1}(-1)$
sub-bundle along which we again double the zero section. We then 
take the ideal sheaf of a reduced point $p\in\PP^1$
and twist
by $\O_{\PP^1}(1)$. 
The resulting sheaf $F$ satisfies $\chi(F)=4$
and has a unique section (up to automorphisms
of the sheaf). We find a $\PP^1\times\PP^1$ of pairs which glues into
$\PP^3$ as the quadric.

However, if the moduli space were really $\PP^3$, then the invariant would
be $-4$, and we would have
 $$Z_{P,2[\PP^1]}(q)=-2q^3-4q^4+\ldots\,$$
in disagreement with the predictions  
$$Z_{GW,2[\PP^1]}(q)= Z'_{DT,2[\PP^1]}(q) =
-2q^3+4q^4+\ldots\,$$ 
of \cite{FaberPand, MNOP1}. In
fact the moduli space has a thickening along the quadric $$\PP^1\times\PP^1
\subset\PP^3$$ corresponding to the first order
movement of the point $p$ 
in the direction of the $\O_{\PP^1}(-1)$ sub-bundle. A computation
shows there
are no more deformations: the moduli space $P_4(X,2[\PP^1])$ 
has Zariski
tangent spaces of dimension 4 along the quadric 
(and 3 over the rest of $\PP^3$).

The natural $(\C^*)^3$ acting on $\O_{\PP^1}(-1)
^{\oplus2}$ fixing the Calabi-Yau form acts on $P_4(X,2[\PP^1])$
 with $4$ isolated fixed points, all lying on the quadric. 
In the virtual localisation formula of \cite{GraberP},
the fixed points each count as $(-1)^4$, where $4$ 
is the dimension of the corresponding
Zariski tangent space. Thus, the extra thickened
dimension along the quadric indeed makes the invariant $4$ instead of $-4$.

\subsection{Contribution of an isolated curve}
Let $X$ be a 3-fold, and let 
$$C\subset X$$
be a nonsingular embedded curve of genus $g$ which represents
an infinitesimally isolated solution of the incidence{\footnote{In
the Calabi-Yau case, no incidence conditions are required.}} conditions
$\prod_i T_i^{e_i}$.
The curve $C$ has a well-defined contribution to the Gromov-Witten 
potential $Z_{GW,[C]}$ by \cite{PandDegen},
\begin{equation}\label{cgj}
Z_{GW,[C]}^C(u) =
\left(\frac{\sin(u/2)}{u/2}\right)^{2g-2+\int_Cc_1(X)}u^{2g-2}.
\end{equation}
We will calculate the contribution of $C$ to the stable pairs theory to
be
\beq{predict}
Z^C_{P,[C]}(q)=q^{1-g}(1+q)^{2g-2+\int_C c_1(X)},
\eeq
in agreement with our conjectures in the Calabi-Yau and Fano cases.

Since all pure rank 1 sheaves on $C$ are locally free,
the stable pairs with support on $C$ are of the form 
$$(F,s) =(\O_C(D),s_D),$$
where $D\subset C$ is a divisor and
 $s_D$ is the canonical section of $\O_C(D)$.  The
cokernel of $s_D$ is $\O_D(D)$.
Thus, the moduli space 
\begin{equation}\label{kkww2}
P^C_{1-g+d}(X,[C])\subset P_{1-g+d}(X,[C])
\end{equation}
of stable pairs cut by the incidence conditions is simply the space 
$\text{Sym}^d(C)$ of degree $d$ divisors on $C$.
The infinitesimal isolation of $C$ implies the inclusion \eqref{kkww2}
is both open and closed.

To calculate the contribution of
$C$, we must identify the obstruction bundle 
$$\text{Obs} \to \text{Sym}^d(C)$$
determined by the deformation theory of the complex 
$$I\udot=\{\O_X\to\O_C(D)\}.$$
In the Calabi-Yau case, 
$$\text{Obs} = T^\vee_{\text{Sym}^d(C)}$$
by the self-duality of the obstruction theory.
Then,
\begin{equation} \label{ZPCq}
Z^C_{P,[C]}(q)\ =\ \sum_{d\geq 0} q^{1-g+d} 
(-1)^d \chi_{top}(\text{Sym}^d(C)).
\end{equation}

\begin{lem}
$\sum_{d\geq 0} \chi_{top}(\text{\em Sym}^d(C))q^d\ =\ (1-q)^{-\chi_{top}(C)}$.
\end{lem}

\begin{proof}
There are many elementary derivations of the result. 
By expressing
$\text{Sym}^d(X)$ as the quotient of $X^{d}$ by the symmetric group,
we see
$\chi_{top}(\text{Sym}^d(X))$ depends only on $\chi_{top}(X)$
for any topological space. 
From the identity
$$
\text{Sym}^d(X\sqcup\{p\})\ \cong\
\text{Sym}^d(X)\ \sqcup\ \text{Sym}^{d-1}(X\sqcup\{p\})\ ,
$$
we deduce 
$$\sum_{d\geq 0} \chi_{top}(\text{Sym}^d(X))q^d\ =
(1-q) \sum_{d\geq 0} \chi_{top}(\text{Sym}^d(X  \sqcup\{p\}     ))q^d\ .$$
Hence, we need only prove the identity for $X$ with positive Euler characteristic,
for which we may replace $X$ by a finite number of points.
\end{proof}

Therefore \eqref{ZPCq} becomes $Z^C_{P,[C]}(q)=q^{1-g}(1+q)^{2g-2}$ and
the verification of \eqref{predict} in the Calabi-Yau
case is complete. 
The stable pairs calculation provides a simple geometric
interpretation of the Hodge integrals entering in \eqref{cgj}.

The obstruction bundle in the  Fano case is determined
in the next result.

\begin{prop} \label{obsbdl}
The obstruction space $\Ext^2(I\udot,I\udot)_0$ sits inside
a canonical exact sequence
$$
0\to H^1(\nu\_C)\to\Ext^2(I\udot,I\udot)_0\to 
H^0(\O_D(D) \otimes K_X)^\vee \to 0
$$
where $\nu\_C$ is the normal bundle to $C\subset X$.
\end{prop}

\begin{proof}
Applying $\Hom(\ \cdot\ ,I\udot)$ to  the exact triangle
$
I\udot\to\O_X\to F
$
yields the exact sequence
\begin{multline} \label{les}
\Ext^2(F,I\udot)\to\Ext^2(\O_X,I\udot)\to\Ext^2(I\udot,I\udot)\to\\
\Ext^3(F,I\udot) \to\Ext^3(\O_X,I\udot)\to\Ext^3(I\udot,I\udot).
\end{multline}
The exact sequence obtained by applying $\Hom(\O_X, \ \cdot \ )$
to the same triangle
and the vanishing $H^2(F)=H^3(F)=0$ together yield
$$\Ext^3(\O_X,I\udot)=\Ext^3(\O_X,\O_X) \cong H^3(\O_X).$$
The last map in \eqref{les} is the identity 
$H^3(\O_X)\to\Ext^3(I\udot,I\udot)$ and an injection
since the composition
$$
H^3(\O_X)\Rt{\id}\Ext^3(I\udot,I\udot)\Rt{\tr}H^3(\O_X)
$$
is multiplication by rank\,($I\udot)=1$. Thus, we may replace \eqref{les}
 by
the exact sequence
\beq{les2}
\Ext^2(F,I\udot)\to\Ext^2(\O_X,I\udot)\to\Ext^2(I\udot,I\udot)\to
\Ext^3(F,I\udot)\to0,
\eeq
fitting into the diagram \vspace{-3mm}
$$
\spreaddiagramcolumns{-.6pc}
\spreaddiagramrows{-1pc}
\xymatrix{
&& 0 \dto \\
\Ext^1(F,F) \rto\dto & H^1(F) \rto\dto & \Ext^2(I\udot,I\udot)_0 \dto \\
\Ext^2(F,I\udot) \rto & \Ext^2(\O_X,I\udot) \rto\dto & \Ext^2(I\udot,I\udot) \rto\dto^{\tr} & \Ext^3(F,I\udot) \rto & 0 \\
& H^2(\O_X) \ar@{=}[r]\dto & H^2(\O_X) \dto \\
& 0 & 0 }
$$
The first two vertical sequences are obtained from 
the exact triangle $I\udot \to \O_X \to  F$
 and  the third is obtained from the
trace map. The diagram commutes because the lower square does.
The commutation for the lower square holds since
 $\Ext^2(\O_X,I\udot) =H^2(\I_C)$ 
and the map to $\Hom(I\udot,I\udot)$ induced by the canonical
maps $\I_C\to I\udot\to\O_X$ takes $f\in\I_C$ to $f\cdot\id$.

The first horizontal map in the above diagram, restricted to 
$$H^1(\O_C)=
H^1(\hom(F,F))\subset\Ext^1(F,F),$$
 is multiplication by the section
$s\colon H^1(\O_C)\to H^1(F)$ which
is onto because $H^1(Q)=0$ since $Q$ is supported in dimension 0.
Therefore the diagram gives
$$
\Ext^2(I\udot,I\udot)_0\ \cong\ \Ext^3(F,I\udot).
$$
The latter
is Serre dual to $\Hom(I\udot,F\otimes K_X)$.  
The exact sequence 
\begin{multline} \label{lirr}
0\to\ext^{-1}(I\udot,F\otimes K_X)\to K_X|_C\rt{s}F\otimes K_X \\
\to\hom(I\udot,F\otimes K_X)\to\ext^1(F,F\otimes K_X)\to0.
\end{multline} is obtained
by applying $\hom(\ \cdot\ ,F\otimes K_X)$
to 
$I\udot\to \O_X \to F$.
Since $s$ is injective, $\ext^i(I\udot,F\otimes K_X)=0$ for
$i<0$ and
\beq{hhom}
\Hom(I\udot,F\otimes K_X)=H^0(\hom(I\udot,F\otimes K_X)).
\eeq
Then, by \eqref{lirr},
$$
0\to\O_D(D)\otimes K_X\to\hom(I\udot,F\otimes K_X)\to\nu\_C\otimes K_X\to0.
$$
Taking $H^0$ we have, by \eqref{hhom},
$$
0\to
H^0(\O_D(D)\otimes K_X)
\to\Hom(I\udot,F\otimes K_X)\to 
H^0(\nu\_C\otimes K_X)\to0.
$$
Dualizing gives the required result.
\end{proof}

Since $C$ is assumed to be infinitesimally isolated, the obstruction
space
$H^1(\nu_C)$ vanishes. Therefore, 
$$\text{Obs}_{[D]} =  H^0(\O_D(D) \otimes K_X)^\vee$$
for $[D]\in \text{Sym}^d(C)$.

Let $Z\subset \text{Sym}^d(C)\times C$ be the universal divisor, and let
$\O_Z(Z)$ be its normal bundle. Let 
$$\pi_d\colon
Z\to \text{Sym}^d(C), \ \ \pi\_C\colon Z\to C$$ 
be the projections. The first projection
 $\pi_d$ is a $d$-fold cover.

Let $L$ be a line bundle on $C$ of degree $l$. A rank $d$ vector bundle
$L_d$ on $\text{Sym}^d(C)$
is obtained tautologically by
$$
L_d=\pi_{d*}(\O_Z(Z) \otimes \pi_C^*L^\vee)^\vee,$$
The identity
$$
\int_{\text{Sym}^d(C)}c_d(L_d)=\left(\!\!\begin{array}{c} 2g-2+l \\ d \end{array}\!\!\right)
=\frac{(2g-2+l)\ldots(2g-2+l-d+1)}{d!}\,,
$$
where we use the formula on the right if $2g-2+l<0$,
is easily obtained via 
intersection theory on the $d!$-fold cover 
$$\epsilon:C^d
\to \text{Sym}^d(C).$$
The product of the $d$ Chern roots of $\epsilon^*(L_d)$ is
$$ \prod_{i=1}^d \left(\pi_i^{*}(c_1(K_C)+c_1(L)) - \sum_{j=i+1}^d
\bigtriangleup_{ij}\right),$$
where 
$$\pi_i: C^d \rightarrow C$$
is the $i^{th}$ projection and $\bigtriangleup_{ij}$ is
the codimension 1 diagonal. 
The leading term 
$$\frac{(2g-2+l)^d}{d!} = \frac{1}{d!}
\int_{C^d} \prod_{i=1}^d \pi_i^{*}\big(c_1(K_C)+c_1(L)\big)$$
is obtained from the leading terms of the Chern roots, and 
the
lower order terms are obtained from the diagonals.

In our case, $\text{Obs}= L_d$ where 
 $L=K_X^\vee|_C$ has degree $l=\int_Cc_1(X)$.  
The stable pairs invariant is
$$
P^C_{1-g+d}(X,[C])=\frac{(2g-2+l)\ldots(2g-2+l-d+1)}{d!}\,.
$$
Therefore
\begin{eqnarray*}
Z_{P,C}(q) &=& \sum_{d\ge0}\frac{(2g-2+l)
\ldots(2g-2+l-d+1)}{d!}q^{1-g+d} \nonumber
\\ &=& q^{1-g}(1+q)^{2g-2+l}\\
& =& q^{1-g}(1+q)^{2g-2+\int_Cc_1(X)},
\end{eqnarray*}
in perfect agreement with (\ref{predict}).

Finally, we see here that the stable pairs invariant is truly 
local with respect to $C$. As in the Gromov-Witten case, the local calculation
\eqref{predict} is valid in any ambient geometry
in which $C$ is infinitesimally isolated. 
The DT
invariant always probes the whole manifold because of the
wandering points. To extract simple local curve results valid
in ambient geometries for DT theory is
technically very difficult \cite{BehrendBryan}.

\section{The vertex} \label{vertex}
\subsection{Overview}
The study of 3-fold theories of nonsingular toric varieties
naturally leads to the notion of a vertex.
The vertex takes its simplest form when restricted to  
the toric Calabi-Yau 
case.{\footnote{Only 
non-compact geometries can be both toric and Calabi-Yau.}} 
In Gromov-Witten theory, the Calabi-Yau vertex
is an evaluation of special Hodge integrals \cite{AKMV}.
The DT vertex, calculated in \cite{MNOP1,MNOP2},
is related to normalized box counting. 
A full development of the stable pairs vertex will be presented
in \cite{PT2}. We give
a short summary of the results and conjectures
in the Calabi-Yau
case here.

\subsection{DT vertex}
Let $\C^3$ have coordinates $x_1,x_2,x_3$. Let the
torus
$$\T=(\C^*)^3$$ acting diagonally on $\C^3$. 
Let $C\subset \C^3$ be a $\T$-fixed 
subscheme of dimension at most 1.
The subscheme $C$ is defined by a 
monomial ideal $$\I_C \subset \C[x_1,x_2,x_3].$$
The latter may be visualized as a 3-dimensional partition $\pi$.
The localisations
$$(\I_C)_{x_1} \subset \C[x_1,x_2,x_3]_{x_1},$$  
$$(\I_C)_{x_2} \subset \C[x_1,x_2,x_3]_{x_2},$$ 
$$(\I_C)_{x_3} \subset \C[x_1,x_2,x_3]_{x_3},$$
are all $\T$-fixed, and each corresponds to a 2-dimensional partition 
$\mu^i$.
Alternatively, the 2-dimensional partitions $\mu^i$ can be defined
as the infinite limits of 
the $x_i$-constant cross-sections of $\pi$.
If all the $\mu^i$ are empty, then $C$ must be 0-dimensional.

Given a triple $\mmuu=(\mu^1,\mu^2,\mu^3)$ of outgoing partitions,
there exists
a unique minimal $\T$-fixed subscheme 
$$C_{\mmuu}\subset \C^3$$
with outgoing partitions $\mu^i$.  
If the $\mu^i$ are not all empty, then
$C_{\mmuu}$ is easily seen to be
a Cohen-Macaulay curve.

Let $S_{\mmuu}$ be the set of $\T$-fixed
subschemes $C\subset \C^3$ with outgoing partitions $\mu^i$.
Since the $\T$-fixed subschemes are isolated, $S_{\mmuu}$
is a discrete set.
Each
$[C]\in S_{\mmuu}$ contains the minimal
subscheme as a quotient
$$\O_C \to {\O_{C_{\mmuu}}} \to 0.$$
The minimal subscheme is the unique Cohen-Macaulay
element of $S_{\mmuu}$.
Let $|\mmuu|$ denote the renormalized volume{\footnote{
The renormalized volume $|\pi|$ is defined by
$$
|\pi| = \#
\left\{\pi \cap [0,\dots,N]^3 \right\}-
(N+1) \sum_1^3 |\mu^i| \,, \quad N\gg 0 \,.
$$
The renormalized volume is
independent of the
cut-off $N$ as long as $N$ is
sufficiently large. The number $|\pi|$
so defined may be negative.
}} 
of the
partition $\pi$ corresponding to $\I_{C_\mmuu}$.
For 
$[C]\in S_{\mmuu}$, define the length by
$$\ell(C)= \dim_\C( {\I_{C_{\mmuu}}}/\I_C) <\infty.$$
Consider the series
$$Z^{DT}_{\mmuu}(q)= 
(-q)^{|\mmuu|}
\sum_{[C]\in S_{\mmuu}} 
(-q)^{\ell(C)}.$$
Since there are only finitely many elements of given length,
$Z^{DT}_{\mmuu}(q)$ is well-defined.

The series $Z^{DT}_{\emptyset,\emptyset,\emptyset}(q)$
is the well-known MacMahon function enumerating finite 3-dimensional
partitions,
$$Z^{DT}_{\emptyset,\emptyset,\emptyset}(q) =\prod_{n\ge1}
\frac{1}{(1-(-q)^n)^{n}}.$$
The normalized vertex \cite{MNOP1,MNOP2}
governing the DT theory in the
toric Calabi-Yau case is
$$\mathsf{W}^{DT}_{\mmuu}(q)= \frac{Z^{DT}_{\mmuu}}
{Z^{DT}_{\emptyset,\emptyset,\emptyset}}(q).$$
The DT vertex can be viewed a normalized count
of $3$-dimensional partitions. 

\subsection{Stable pairs vertex}
Given a $\T$-fixed stable pair on $\C^3$,
$$\O_{\C^3} \stackrel{s}{\to} F,$$ 
the scheme theoretic support must be of the form
$$C_{\mmuu}\subset \C^3.$$
Moreover,
the quotient $Q$ must be supported at the origin.

We will use the characterization of
 Proposition \ref{desc} to study $\T$-fixed stable pairs with
support $C_{\mmuu}$. 
Let
$$\m\subset \O_{C_{\mmuu}}$$
be the ideal sheaf of the origin in $C_\mmuu$.
Let 
$$M_1 = \C[x_1,x_2,x_3]_{x_1}/\I_{\mu^1},$$  
$$M_2=   \C[x_1,x_2,x_3]_{x_2}/\I_{\mu^2},$$ 
$$M_3= \C[x_1,x_2,x_3]_{x_3}/\I_{\mu^3},$$
be the quotients by the 
 $\T$-fixed ideals $\I_{\mu^i}$ determined by the respective outgoing
partitions, and let
$$M_{\mmuu}= \oplus_{i=1}^3 M_i$$
be viewed as a $\C[x_1,x_2,x_3]$-module.
There is a canonical homomorphism
$$\O_{C_{\mmuu}} \to M_{\mmuu}$$
given by 
$$ 1 \mapsto (1,1,1).$$ 
The limit in Proposition \ref{desc} has a simple
description,
$$\lim_{\To}\hom(\m^r,\O_{C_{\mmuu}})/\O_{C_\mmuu} = M_\mmuu/\O_{C_\mmuu}$$
as a $\C[x_1,x_2,x_3]$-module.

By Proposition \ref{desc}, $\T$-fixed stable pairs on $\C^3$
correspond to 
$\T$-fixed
$\C[x_1,x_2,x_3]$-submodules $Q\subset M_\mmuu/\O_{C_\mmuu}$ of finite length
$$\ell(Q) = \dim_\C(Q).$$
Unlike the DT case, such $\T$-fixed submodules 
{\em need not be isolated}. However, the geometry of the moduli of
$\T$-fixed submodules modules of
$M_\mmuu/\O_{C_\mmuu}$
 is elementary: the components are simply 
products{\footnote{The $0^{th}$ product is a point. In fact,
the $\T$-fixed modules are isolated if at least one of the $\mu^i$ is
empty. Indeed, the proof Conjecture \ref{kqqwe} in the
1 and 2-leg case is not difficult \cite{PT2}. 
Non-trivial moduli of $\T$-fixed submodules
 appear only in the full $3$-leg case.}} of $\PP^1$.

Let the set 
$S^M_\mmuu$ index the components $\mathcal{Q}$ of the
moduli space  of $\T$-fixed
submodules of $M_\mmuu/\O_{C_\mmuu}$.
The components are easily indexed by a
box counting strategy explained in \cite{PT2}.
The stable pairs vertex is conjectured{\footnote{Proofs
are given in the 1 and 2-leg cases.}} in \cite{PT2} to be
the series
$$\mathsf{W}^{P}_{\mmuu}(q)= 
(-q)^{|\mmuu|}
\sum_{[\mathcal{Q}]\in S^M_{\mmuu}} 
\chi_{top}(\mathcal{Q})(-q)^{\ell(\mathcal{Q})}.$$
The topological Euler characteristic above is always
a power of 2.
Since the length is constant in components, $\ell(\mathcal{Q})$
is well-defined.
Again,
there are only finitely many components of given length.

We conjecture a correspondence between the DT
and stable pairs vertices.

\begin{conj}\label{kqqwe}
$\mathsf{W}^{DT}_{\mu^1,\mu^2,\mu^3}(q) = 
\mathsf{W}^{P}_{\mu^1,\mu^2,\mu^3}(q).$
\end{conj}

Independent of the geometric and string theoretic framework,
Conjecture \ref{kqqwe} is a striking statement about the
algebro-combinatorics of $\C[x_1,x_2,x_3]$.

\subsection{Example}
The simplest 3-leg 
example occurs when $$\mu^1=\mu^2=\mu^3=(1).$$
 The 
 Cohen-Macaulay curve $C_{(1),(1),(1)}$ 
is the union
of the 3 coordinate axes.  From the definitions, we find
\begin{equation}
\label{jjj}
\frac{M_{(1),(1),(1)}}{\O_{C_{(1),(1),(1)}}} =
\frac{\C[x_1,x_1^{-1}] \oplus \C[x_2,x_2^{-1}] \oplus \C[x_3,x_3^{-1}]}
{(1,1,1)\cdot \C[x_1,x_2,x_3]}.
\end{equation}

To compute the stable pairs vertex, we must study the
$\T$-fixed submodules of \eqref{jjj}. Of course,  
 0 is the unique submodule of length 0.
The first interesting case is length 1.
Then, 
$$(\gamma_1,\gamma_2,\gamma_3) \in \ \C^3 \setminus \langle \ (1,1,1)
\ \rangle$$
determines a $\T$-fixed submodule of length 1 of \eqref{jjj}.
So the moduli space for length 1 is isomorphic to $\PP^1$.
For length 2, the $\T$-fixed submodules are
$$\langle\ (x^{-1}_1,0,0)\ \rangle, \ \ \langle\ (0,x_2^{-1},0)\ \rangle, 
\ \ \langle \ (0,0,x_3^{-1}) \ \rangle,$$ 
$$\langle \ (1,0,0),\  (0,1,0),\ (0,0,1)\ \rangle$$
all of which are isolated points.
In general, the $\T$-fixed submodules are elementary to enumerate
since the generators must be vectors of monomials.
We find
\begin{eqnarray*}
(-q)^{2}\cdot \mathsf{W}^P_{(1),(1),(1)}(q) & =& 1 + 2(-q) + \sum_{n\geq 2} \left(
\binom{n}{2}+3\right) (-q)^n \\
& = & \frac{1-q^5}{(1-q)(1+q)^3} \\
& = & 1 -2q +4 q^4 -6 q^3 + \ldots
\end{eqnarray*}
Indeed, the vertex $\mathsf{W}^P_{\mmuu}(q)$ is always
a rational function.

The DT vertex $\mathsf{W}^{DT}_{(1),(1),(1)}$ counts
boxes added to the Cohen-Macaulay curve $C_{(1),(1),(1)}$.
We find
$$(-q)^2 \cdot \mathsf{W}^{DT}_{(1),(1),(1)}(q)  = 
1-3q+9q^2 -22 q^3 + \ldots.$$
Since the MacMahon function is
$$\mathsf{W}^{DT}_{\emptyset,\emptyset,\emptyset}(q) =
1-q + 3q^2 - 6 q^3 +\ldots,$$
Conjecture \ref{kqqwe} is verified to order 3 by
$$1 -2q +4 q^4 -6 q^3 + \ldots =
\frac{1-3q+9q^2 -22 q^3 + \ldots}
{1-q + 3q^2 - 6 q^3+ \ldots} \, .$$
In fact, the exact equality of Conjecture \ref{kqqwe} is
not hard to obtain by box counting for the 
example $((1),(1),(1))$.

\vspace{+8 pt}
\noindent
Department of Mathematics\\
Princeton University\\
rahulp@math.princeton.edu

\vspace{+8 pt}
\noindent
Department of Mathematics \\
Imperial College \\
rpwt@imperial.ac.uk


\begin{thebibliography}{MNOP2}

\bibitem{AKBranch}
V.~Alexeev and A.~Knutson.
\newblock {\em Complete moduli spaces of branchvarieties}, 
\newblock math.AG/0602626.

\bibitem{AKMV}
M.~Aganagic, A.~Klemm, M.~Mari{\~n}o, and C.~Vafa.
\newblock {\em The topological vertex}, Comm. Math. Phys., {\bf 254}, 425--478,
  2005.
\newblock hep-th/0305132.

\bibitem{AM}
P.~S. Aspinwall and D.~R. Morrison.
\newblock {\em Topological field theory and rational curves}, Comm. Math.
  Phys., {\bf 151}, 245--262, 1993.

\bibitem{BayerPoly} A.~Bayer, \newblock{\em Polynomial {B}ridgeland stability
conditions and the large volume limit}, arXiv:0712.1083.

\bibitem{BehrendDT}
K.~Behrend.
\newblock {\em {D}onaldson-{T}homas invariants via microlocal geometry},
to appear in Ann. of Math., 2009.
\newblock math.AG/0507523.

\bibitem{BehrendBryan}
K.~Behrend and J.~Bryan.
\newblock {\em Super-rigid {D}onaldson-{T}homas invariants},
Math. Res. Lett.  {\bf 14}, 559--571, 2007.
\newblock math.AG/0601203.

\bibitem{BFNormalCone}
K.~Behrend and B.~Fantechi.
\newblock {\em The intrinsic normal cone}, Invent. Math., {\bf 128}, 45--88,
  1997.
\newblock math.AG/9601010.

\bibitem{BFHilb}
K.~Behrend and B.~Fantechi.
\newblock {\em Symmetric obstruction theories and {H}ilbert schemes of points
  on threefolds}, Algebra Number Theory {\bf 2}, 313--345, 2008.
\newblock math.AG/0512556.

\bibitem{BryanP}
J.~Bryan and R.~Pandharipande.
\newblock {\em The local {G}romov-{W}itten theory of curves}, J.
  Amer. Math. Soc.  {\bf 21}, 101--136, 2008.
\newblock math.AG/0411037.

\bibitem{BridgelandFlops}
T.~Bridgeland.
\newblock {\em Flops and derived categories}, Invent. Math., {\bf 147},
  613--632, 2002.
\newblock math.AG/0009053.

\bibitem{BrStability}
T.~Bridgeland.
\newblock {\em Stability conditions on triangulated categories}, 
  Ann. of Math. {\bf 166}, 317--345, 2007.
\newblock math.AG/0212237.

\bibitem{COGP}
P.~Candelas, X.~C. de~la Ossa, P.~S. Green, and L.~Parkes.
\newblock {\em A pair of {C}alabi-{Y}au manifolds as an exactly soluble
  superconformal theory}, Nuclear Phys. B, {\bf 359}, 21--74, 1991.


\bibitem{DenefMoore}
F.~Denef and G.~W. Moore.
\newblock {\em Split states, entropy enigmas, holes and halos},
\newblock hep-th/0702146.

\bibitem{DiacADHM}
D.~E. Diaconescu.
\newblock {\em Moduli of {ADHM} sheaves and local {D}onaldson-{T}homas theory},
  arXiv:0801.0820.



\bibitem{DiacMoore}
D.~E. Diaconescu and G.~W. Moore.
\newblock {\em Crossing the wall: Branes vs. bundles},
\newblock hep-th/07063193.

\bibitem{DT}
S.~K. Donaldson and R.~P. Thomas.
\newblock {\em Gauge theory in higher dimensions}.
\newblock In {\em The geometric universe (Oxford, 1996)},  31--47. Oxford Univ.
  Press, Oxford, 1998.

\bibitem{FaberPand}
C.~Faber and R.~Pandharipande.
\newblock {\em Hodge integrals and {G}romov-{W}itten theory}, Invent. Math.,
  {\bf 139}, 173--199, 2000.
\newblock math.AG/9810173.

\bibitem{GraberP}
T.~Graber and R.~Pandharipande.
\newblock {\em Localization of virtual classes}, Invent. Math., {\bf 135},
  487--518, 1999.
\newblock math.AG/9708001.

\bibitem{GV1}
R.~Gopakumar and C.~Vafa.
\newblock {\em M-theory and topological strings--{I}}, 
\newblock hep-th/9809187.

\bibitem{GV2}
R.~Gopakumar and C.~Vafa.
\newblock {\em M-theory and topological strings--{II}}, 
\newblock hep-th/9812127.


\bibitem{Honsen}
M.~Honsen.
\newblock {\em A compact moduli space parameterizing Cohen-Macaulay curves
in projective space}.
\newblock PhD thesis, MIT, 2004.


\bibitem{hetal}
S.~Hosono, M.-H. Saito, and A. Takahashi.
\newblock {\em Relative Lefschetz actions and BPS state
counting}, IMRN, {\bf 15}, 783--816, 2001.

\bibitem{HLShaves}
D.~Huybrechts and M.~Lehn.
\newblock {\em The geometry of moduli spaces of shaves}.
\newblock Aspects of Mathematics, E31. Friedr. Vieweg \& Sohn, Braunschweig,
  1997.


\bibitem{HT}
D.~Huybrechts and R.~P.~Thomas.
\newblock {\em Deformation-obstruction
theory for complexes via Atiyah and Kodaira--Spencer classes},
arXiv:0805.3527.

\bibitem{InabaModuli}
M.~Inaba.
\newblock {\em Moduli of stable objects in a triangulated category},
\newblock math.AG/0612078.

\bibitem{IonelParker}
E.~Ionel and T.~H. Parker.
\newblock {\em Relative {G}romov-{W}itten invariants}, Ann. of Math. (2), {\bf
  157}, 45--96, 2003.
\newblock math.SG/9907155.

\bibitem{JoWallCrossing}
D.~Joyce.
\newblock {\em Configurations in abelian categories. {IV}. {I}nvariants and
  changing stability conditions}, Adv. Math. {\bf 217}, 125--204, 2008.
\newblock math.AG/0410268.

\bibitem{KatzBPS}
S.~Katz.
\newblock {\em Genus zero {G}opakumar-{V}afa invariants of contractible
  curves},  J. Differential Geom. {\bf 79}, 185--195, 2008.
\newblock math.AG/0601193.

\bibitem{KKVSpinning}
S.~Katz, A.~Klemm, and C.~Vafa.
\newblock {\em M-theory, topological strings and spinning black holes}, Adv.
  Theor. Math. Phys., {\bf 3}, 1445--1537, 1999.
\newblock hep-th/9910181.


\bibitem{KlemmP}
A.~Klemm and R.~Pandharipande.
\newblock {\em Enumerative geometry of {C}alabi-{Y}au 4-folds}, Comm. Math.
Phys. {\bf 281}, 621--653, 2008.
\newblock math.AG/0702189.


\bibitem{KollarJDG}
J.~Koll{\'a}r.
\newblock {\em Projectivity of complete moduli}, J. Differential Geom., {\bf
  32}, 235--268, 1990.

\bibitem{LPPairs1}
J.~Le~Potier.
\newblock {\em Syst\`emes coh\'erents et structures de niveau}, Ast\'erisque,
  {\bf 214}, 143, 1993.

\bibitem{LPPairs2}
J.~Le~Potier.
\newblock {\em Faisceaux semi-stables et syst\`emes coh\'erents}.
\newblock In {\em Vector bundles in algebraic geometry (Durham, 1993)}, vol
208
  of London Math. Soc. Lecture Note Ser.,  179--239. Cambridge Univ. Press,
  Cambridge, 1995.

\bibitem{LPCobordism}
M.~Levine and R.~Pandharipande.
\newblock {\em Algebraic cobordism revisited}, 
\newblock math.AG/0605196.

\bibitem{LiRuan} 
A.-M.~Li and Y.~Ruan. 
\newblock \emph{Symplectic surgery and
Gromov-Witten invariants of Calabi-Yau 3-folds I},  
Invent.\ Math.\ \textbf{145}, 151--218, 2001.

\bibitem{LiRelative}
J.~Li.
\newblock {\em Stable morphisms to singular schemes and relative stable
  morphisms}, J. Differential Geom., {\bf 57}, 509--578, 2001.
\newblock math.AG/0009097.

\bibitem{LiRelativeGW}
J.~Li.
\newblock {\em A degeneration formula of {GW}-invariants}, J. Differential
  Geom., {\bf 60}, 199--293, 2002.
\newblock math.AG/0110113.

\bibitem{LiDT}
J.~Li.
\newblock {\em Zero dimensional {D}onaldson-{T}homas invariants of threefolds},
 Geom. Topol. {\bf 10}, 2117--2171, 2006. \newblock math.AG/0604490.

\bibitem{LiTianVirtual}
J.~Li and G.~Tian.
\newblock {\em Virtual moduli cycles and {G}romov-{W}itten invariants of
  algebraic varieties}, J. Amer. Math. Soc., {\bf 11}, 119--174, 1998.
\newblock math.AG/9602007.

\bibitem{LiWuDegenDT}
J.~Li and B.~Wu.
\newblock {\em Degeneration of {D}onaldson-{T}homas invariants}, in
  preparation.

\bibitem{LieblichModuli}
M.~Lieblich.
\newblock {\em Moduli of complexes on a proper morphism}, J. Algebraic Geom.,
  {\bf 15}, 175--206, 2006.
\newblock math.AG/0502198.

\bibitem{Lowen}
W.~Lowen.
\newblock {\em Obstruction theory for objects in abelian and derived
categories}, Comm.\ Algebra {\bf 33} (2005), 3195--3223.
\newblock math.KT/0407019.

\bibitem{MNOP1}
D.~Maulik, N.~Nekrasov, A.~Okounkov, and R.~Pandharipande.
\newblock {\em Gromov-{W}itten theory and {D}onaldson-{T}homas theory. {I}},
  Compos. Math., {\bf 142}, 1263--1285, 2006.
\newblock math.AG/0312059.

\bibitem{MNOP2}
D.~Maulik, N.~Nekrasov, A.~Okounkov, and R.~Pandharipande.
\newblock {\em Gromov-{W}itten theory and {D}onaldson-{T}homas theory. {II}},
  Compos. Math., {\bf 142}, 1286--1304, 2006.
\newblock math.AG/0406092.

\bibitem{PandDegen}
R.~Pandharipande.
\newblock {\em Hodge integrals and degenerate contributions}, Comm. Math.
  Phys., {\bf 208}, 489--506, 1999.
\newblock math.AG/9811140.

\bibitem{PandICM}
R.~Pandharipande.
\newblock Three questions in {G}romov-{W}itten theory.
\newblock In {\em Proceedings of the International Congress of Mathematicians,
  Vol. II (Beijing, 2002)} 503--512, Beijing, 2002. Higher Ed. Press.
  math.AG/0302077.

\bibitem{PT2}
R.~Pandharipande and R.~P. Thomas.
\newblock {\em The 3-fold vertex via stable pairs}, to appear in Geom. Topol. arXiv:0709.3823.

\bibitem{PZ}
R.~Pandharipande and A.~Zinger. {\em Enumerative geometry of
Calabi-Yau 5-folds}, arXiv:0802.1640.


\bibitem{SzendroiDT}
B.~Szendr{\H o}i.
\newblock {\em Non-commutative {D}onaldson-{T}homas theory and the conifold},
 Geom. Topol. {\bf 12}, 1171--1202, 2008. 
 \newblock math.AG/07053419.

\bibitem{ThCasson}
R.~P. Thomas.
\newblock {\em A holomorphic {C}asson invariant for {C}alabi-{Y}au 3-folds, and
  bundles on {$K3$} fibrations}, J. Differential Geom., {\bf 54}, 367--438,
  2000.
\newblock math.AG/9806111.

\bibitem{Toda}
Y.~Toda.
\newblock {\em Birational {C}alabi-{Y}au 3-folds and
BPS state counting}, Commun. Number Theory Phys. {\bf 2}, 63--112, 2008.
\newblock math.AG/07071643.

\bibitem{TodaStab}
Y.~Toda,
\newblock {\em Limit stable objects on Calabi-Yau 3-folds,}
\newblock arXiv:0803.2356.

\end{thebibliography}
\end{document}